\documentclass[letterpaper, 10pt]{article}

\pdfoutput=1

\usepackage{amsmath,amssymb,amsthm}
\usepackage{enumerate}
\usepackage{graphicx}
\usepackage{color}
\usepackage{vmargin}

\input{perso.sty}

\title{Multigrid methods for two-player zero-sum stochastic games}

\author{Marianne Akian \and Sylvie Detournay}

\begin{document}

\maketitle

\begin{abstract}
We present a fast numerical algorithm for large scale zero-sum stochastic 
games with perfect information, which combines policy iteration and 
algebraic multigrid methods. This algorithm can be applied either 
to a true finite state space 
zero-sum two player game or to the discretization of an Isaacs equation.
We present numerical tests 
on discretizations of Isaacs equations or variational inequalities. 
We also present a full multi-level policy iteration, 
similar to FMG, which allows to improve substantially
the computation time for solving some variational inequalities.
\end{abstract}

\section{Introduction}

In the present paper, we are interested in solving non-linear finite
dimensional equations of the form~:
\begin{equation} \label{eq1intro}
v(x) \, = \, {\max_{a \in \A} \,\left(  \min_{b \in \B} \, \left( \sum_{y \in \X}
\mu \, p(y\, | \, x, a, b) \, v(y) \, + \, r(x,a, b) \right) \right)} \qquad \forall x \in \X.
\end{equation}
with unknown the function $v : \X \rightarrow \R$, 
where $\X := \set{1, \dots, n}$.
Here $\A : = \set{1, \dots, m_1}$, $\B : =
\set{1, \dots, m_2}$ are finite sets; 
the functions $(x,a,b) \in \X \times \A \times \B \rightarrow r(x,a,b) \in \R$ and
$(x,a,b,y) \in \X \times \A \times \B \times \X \rightarrow p(y|x,a,b)
\in \R_+$ are given such that $\sum_{y\in\X} p(y|x,a,b) =
1$; and $0<\mu<1$ is a given constant. 
These equations appear when solving the following particular dynamic games.

An infinitely repeated game, or discrete time dynamic game, or
infinite horizon multi-stage game, consists in an infinite sequence of
state transitions, where at each step,  the transition depends on the
actions of the players, and each player receives a reward 
which depends on the state of the game and the actions of all players
at this step.
The aim of each player is to maximize his own objective function, for
instance his payoff which is the sum of the rewards he received at all steps.
The game is stochastic when the state sequence is a random process
with a Markov property, then the objective function is the expected
payoff.
It is a two player zero-sum game when there are two players with
opposite rewards, hence player $2$ aims to minimize
player $1$ objective function.
When the game does not stop in finite time (almost surely), one often
consider a discounted payoff where the reward at each step $k$ is
discounted by some multiplicative factor $\mu^k$, with $0<\mu<1$.

Consider in particular a two player zero-sum discounted stochastic game with
finite state space $\X$ and action spaces $\A$ and
$\B$ for player $1$ and player $2$ respectively.
Denote by $r(x,a,b)$ the reward of player $1$  when (at the current
step) the state is $x\in\X$ and the 
actions of player $1$ and $2$ are $a\in\A$, $b\in\B$ respectively.
Denote by $p(y|x,a,b)$ the transition probability from state $x$
to state $y$  when the actions of player $1$ and $2$ are
$a\in\A$, $b\in\B$ respectively.
Assume that player $1$ plays before player $2$, and that at each step, 
player $1$ is choosing his action $a\in\A$ as a function of the
current state $x\in \X$, and player $2$ is choosing his action
$b\in\B$ as a function of the current state $x\in\X$ and action
$a\in\A$ of player $1$.
Assume each player is maximizing his own objective function. 
Under the previous finiteness conditions on ($\X$, $\A$, $\B$), there 
exists a function $v : \X \rightarrow \R$ which associates to each
$x\in\X $ the expected payoff $v(x)$ of player $1$ when the initial
state of the game is $x$. This function is called the value or the value function
of the game. It is the unique solution of
Equation~\eqref{eq1intro} \cite{Shapley}, called itself
the dynamic programming equation or Shapley equation of the
game. 
Solving Equation~\eqref{eq1intro} is important since it also gives the 
optimal stationary strategies of the game, see section~\ref{discrete} for precise
definitions of strategies and details.  
Discrete time zero-sums stochastic games arise in several domains of
applications, such as military operations~\cite{mceneaneyaircraft04}, 
network flow control~\cite{altmanflow94}, pursuit-evasion problems
(although often studied in the deterministic case), 
see \cite{SorinNeumann} for other applications references.

Equations of the form~\eqref{eq1intro} can also be obtained as
special discretizations of partial differential equations associated to
differential stochastic games, where the state space $\X$ is now a subset of
$\R^d$ (see section~\ref{section-SDG} for details).
For instance the following non-linear elliptic partial differential
equation called Isaacs equation~:
\begin{equation} \label{eq2intro}
\max_{a \in \A} \, \left(  \min_{b \in \B} \, \left( 
\sum_{i,j  = 1}^d \ai_{ij}(a, b, x) \hess{v(x)}{x_i}{x_j}
 +  \sum_{j=1}^d \bi_{j}(a,b,x) \der{v(x)}{x_{j}}  - \lambda v(x)
 + r(x,a,b) \,\right) \right) = 0 \;\;\forall x \in \X
\end{equation}
allows one to solve a differential game in the same way as~\eqref{eq1intro}
solves a discrete time dynamic game.
Here $\A$, $\B$ are either finite sets or subsets of some $\R^p$
spaces, $\lambda \ge 0$ is a scalar, and
$(x , a, b) \in \X \times \A \times\B \rightarrow \ai(x , a, b) = (\ai_{ij}(x, a,
b))_{i,j =1,\dots,d}$ $\in S^+_d$, the set of positive definite
symmetric $d\times d$ matrices, $(x ,a, b) \in \X \times \A \times\B
\rightarrow \bi(x , a, b) = (\bi_{j}(x,a,b))_{j=1,\dots,d}$ $\in
\R^d$, and $(x ,a, b) 
\in \X \times \A \times\B \rightarrow r(x , a, b) \in \R$
are given functions.
Such equations may be applied in particular to pursuit-evasion games 
(see for instance \cite{BaFaSo94}), but they also appear
in solving $H^\infty$ optimal control problems
(see for instance~\cite{basarbernhard}),
or risk-sensitive optimal control problems~\cite{fleming06},
in particular for finance applications~\cite{elliott11}.
The discretization of Equation~\eqref{eq2intro} with a 
monotone scheme in the sense of~\cite{barlsoug} yields an 
equation of the form~\eqref{eq1intro} which can then 
be interpreted as the dynamic 
programming equation of a stochastic game with discrete time and 
finite state space. Suitable possible discretizations schemes are for
instance~: Markov chain discretizations~\cite{kus77,KuDu92}, monotone 
discretizations~\cite{barlsoug}, full discretizations of
semi-Lagrangian type~\cite{BaFaSo94}, and  max-plus finite element
method~\cite{akgbla08} for deterministic games or control problems.
Hence, we are interested in solving discretizations of
Equation~\eqref{eq2intro} which have the form of 
Equation~\eqref{eq1intro}, in order to find an 
approximation of the value of the corresponding differential stochastic game.

In the presence of a discount factor $\mu < 1$, the nonlinear
equation~\eqref{eq1intro} can be solved by applying the fixed point
iterations which are called, in the optimal control and game literature,
value iterations or the value iteration algorithm~\cite{Bellman}. 
The iterations of this method are cheap but their convergence
slows considerably  as the discount factor $\mu$ approaches one. 
Moreover, when we discretize Equation~\eqref{eq2intro} with a finite
difference or finite element method with a discretization step $h$,
we obtain an equation of the form~\eqref{eq1intro} with a
discount factor 
$\mu=1-O(\lambda h^2)$, then when $h$ is small $\mu$ 
is close to one
and the value iteration method is as slow as the Jacobi or
Gauss-Seidel iterations for a discretized linear elliptic equation.
Another approach consists in the so called policy 
iteration algorithm, initially introduced by Howard~\cite{Howard60}
for one player stochastic games (i.e.\ stochastic control problems). Later
adaptations of this algorithm were proposed for the two player games~: by Hoffman and
Karp~\cite{HoffmanKarp} for a special mean-payoff case, by
Dernado~\cite{Denardo} for approximations of value functions in
discounted stochastic games, in Puri thesis~\cite{Puri} for discounted
stochastic games, and by 
Cochet-Terrasson and Gaubert~\cite{CochetGaub} for the general
mean-payoff case. 
In all cases, policy algorithm converges faster than the
value iteration algorithm and in practice it ends in few steps (see
for instance~\cite{DhinGaub} for numerical examples in the case
of deterministic games).

A (feedback) policy (or pure Markovian stationary strategy, see 
Section~\ref{discrete} below)
${\balpha}: \X\rightarrow \A$ for the first player is a 
function which maps any $x\in \X$ to an action $a \in
\A$. Then, starting with an initial policy for player $1$, the policy
iteration algorithm for the two player zero-sum stochastic game consists
in applying successively a policy evaluation step followed by
a policy improvement step. The policy evaluation step amounts to
compute the value of the game for the current policy $\balpha$,
that is the solution $v$ of~\eqref{eq1intro} where instead of taking 
the maximum of the expression inside the ``max'', one evaluates it 
with $a=\balpha(x)$. The policy improvement step consists in finding 
the optimal policy for the current value function $v$,
that is the policy optimizing the expression inside the ``max'' 
in~\eqref{eq1intro} when the value function is $v$.
Computing the above value functions (in the policy evaluation steps) is
performed using the policy iteration algorithm for a one-player game.
The policy iteration algorithm is explained in more general settings
in Section~\ref{section-backG}.  
It stops after a finite number of steps
when the sets of actions are finite,
see~\cite{LionsMercier80,Bertsekas87,puterman} for one player games
and~\cite{Puri,CochetGaub} for two player games.
In addition,  under regularity assumptions on the maps $r$ and $p$, the policy 
iteration algorithm for a one player game with infinite action spaces
is equivalent to Newton's method, thus can have a
super-linear convergence in the neighborhood of the solution,
see~\cite{PutBrum,zidaniBokanowski09} for superlinear convergence under
general regularity assumptions, and~\cite{PutBrum,aki90b,BankRose82}
for order $p>0$ superlinear convergence 
under additional regularity and strong convexity assumptions.



Each policy iteration for a one player game (or each iteration in the
inner loop of the two player algorithm) requires the solution of a
linear system. Indeed, when we fix feedback policies $\balpha : \X
\rightarrow  \A$ and $\bbeta : \X \rightarrow \B$ for player $1$
and $2$ respectively,
the system of equations~\eqref{eq1intro} yields a linear system of
the form~: $v = \mu M v + r$ where $v, r \in \R^\X$ are respectively
the value function of the game and the vector of rewards for the fixed
policies ${\balpha}$ and $\bbeta$,
$0<\mu < 1$ is the discount factor and $M \in \R^{\X \times \X}$ is a Markov matrix
whose elements are the transition probabilities ${M}_{xy} =
p(x|y,\balpha(x),\bbeta(x)) \in \R_+$ for $x,y\in\X$ (and each rowsum of $M$ equals one). 
When the dynamic programming equation~\eqref{eq1intro} is coming from
the discretization of an Isaacs partial differential
equation~\eqref{eq2intro}, this linear system
corresponds to the discretization of a linear elliptic partial differential
equation, hence it may be solved in the best case in a time in the 
number of discretization points  by using multigrid methods, 
that is  the cardinality $|\X|$ of the 
discretized state space $\X$, or the size of the matrix $M$.
For general stochastic games on a finite state space $\X$, since $M$ is a
Markov matrix, the matrix $(I - \mu M)$ of the
linear system is an invertible M-matrix~\cite{plemmons}, and one may
expect the same complexity when solving them by using an algebraic
multigrid method.


In the present paper, we consider the combination of policy iterations 
with the algebraic multigrid method (AMG) introduced by
Brandt, McCormick and Ruge~\cite{Brandt1,Brandt}, see also Ruge and
St\"uben~\cite{stub1}. We shall call AMG$\pi$ the resulting algorithm.
This algorithm can be applied either to a true finite state space 
zero-sum two player game or to the discretization of an Isaacs equation, 
although in the present paper we restrict ourselves to numerical tests
for the discretization of stochastic differential games, since the AMG
algorithm needs some improvements to be applied to arbitrary non
symmetric linear systems arising in game problems.  
Such an association of multigrid methods with policy iteration has
already been used and studied in the case of one player games, that is
discounted stochastic control problems (see Hoppe~\cite{Hope86,Hope87}
and Akian~\cite{akian,aki90b} for Hamilton-Jacobi-Bellman equations
or variational inequalities, Ziv and Shimkin~\cite{Ziv} for AMG with
learning methods). However, it is new in the case of two player
games. 
We have implemented this algorithm (in C) and shall present numerical
tests on discretizations of Isaacs or Hamilton-Jacobi-Bellman
equations or variational inequalities, while comparing AMG$\pi$ with
the combination of policy iterations with direct solvers.

The complexity of two player zero-sum stochastic games is
still unsettled, one only knows that it belongs to the complexity
class of NP$\cap$coNP~\cite{Puri}. 
Indeed, the number of policy iterations is bounded by the 
number of possible policies, which is exponential in the cardinality
of $\X$. Friedmann has shown~\cite{Friedmann} that a strategy improvement
algorithm requires an exponential number of iterations for a ``worst''-case
family of games called parity games, this result can be extended to
other types of zero-sum stochastic games, in particular to mean-payoff and discounted zero-sum stochastic games, and to undiscounted stochastic control problems 
(one-player games) as shown by Fearnley~\cite{fearnley2,fearnley}. 
However, as for Newton's algorithm, convergence can be improved by
starting the policy iteration with a good initial guess, close to the
solution. With this in mind, we present a full multi-level policy
iteration, similar to FMG.  
It consists in solving the problem at each grid level
by performing policy iterations 
until a convergence criterion is verified, 
then to interpolate the strategies and 
value to the next level, in order to initialize the policy iterations 
of the next level, until the finest level is attained. 
When at each level 
policy iterations are combined with the algebraic multigrid method,
we shall call FAMG$\pi$ the resulting full multi-level policy
iteration algorithm.
For one-player discounted games with infinite number of actions and
under regularity assumptions, one can show~\cite{aki90b,akian}
that this kind of full multi-level policy iteration has a computing time 
in the order of the cardinality $|\X|$ of the discretized state space
$\X$ at the finest level.
In Section~\ref{section-numerics}, we give numerical examples on variational inequalities for
two player games, the computation time of which is 
improved substantially using FAMG$\pi$ instead of AMG$\pi$.


The paper is organized as follow. The three following sections are
some recalls about basic definitions on the subject. In Section~\ref{discrete},
we introduce the definition of a two player zero-sum stochastic game
with finite state space and the corresponding dynamic programming
equation. Section~\ref{section-SDG} is about two player zero-sum stochastic
differential games, we recall here the definition of the Isaacs
equation, the variational inequalities and the discretization scheme
that we use. Section~\ref{section-backG} is devoted to the numerical background needed
to solve the dynamic programming equation, including the policy
iteration algorithm and the algebraic multigrid method. Section~\ref{section-AMGpi}
describes our algorithms AMG$\pi$ and FAMG$\pi$. We
present in Section~\ref{section-numerics} some numerical tests on 
discretizations of Isaacs equations and variational inequalities. 
Last section gives ending remarks.


\section{Two player zero-sum stochastic games: the discrete case}
\label{discrete}

The class of two player zero-sum stochastic game was first introduced
by Shapley in the early fifties \cite{Shapley}. We recall in this
section the definition of these games in the case of finite state
space and discrete time (for more details see \cite{Shapley,FilarVrieze,Sorin}). 


We consider a finite state space $\X = \{1, \dots, n\}$. 
A  stochastic process $\proc{\Xk_k}$ on $\X$ gives the state of the
game at each point time $k$, called stage. At each of these stages,
both players have the possibility to influence the course of the
game. 

The stochastic game $\Gamma(x_0)$ starting from $x_0 \in \X$ is played
in stages as follows. The initial state $\Xk_0$ is equal to $x_0$ and known by
the players. The player who plays first, say \MAX, chooses an action
$\Ak_0$ in a set of possible actions $\A(\Xk_0)$. Then the second player,
called \MIN\, chooses an action $\Bk_0$ in a set of possible actions
$\B(\Xk_0,\Ak_0)$. The actions of both players and the current state
determine the payment $r(\Xk_0,\Ak_0,\Bk_0)$ made by \MIN\ to \MAX\ and the
probability distribution $p(\cdot | \Xk_0,\Ak_0,\Bk_0)$ of the new state
$\Xk_1$. Then the game continues in the same way with state $\Xk_1$ and so
on.

At a stage $k$, each player chooses an action knowing the history
defined by $\Ik_k = (\Xk_0, \Ak_0, \Bk_0, \cdots, \Xk_{k-1}, \Ak_{k-1}, \Bk_{k-1},
\Xk_k)$ for \MAX\ and $(\Ik_k, \Xk_k)$ for \MIN. We call a strategy or
policy for a player, a rule which tells him the action to choose 
at any stage and in any situation. 
There are several classes of strategies. 
Assume $\A(x) \subset \A$ and $\B(x,a) \subset \B$ for some sets $\A$ and $\B$.
A behavior or 
randomized strategy for \MAX\ (resp.\ \MIN) is a
sequence $\salpha := (\falpha_0, \falpha_1, \cdots)$ (resp.\ $\sbeta :=
(\fbeta_0, \fbeta_1, \cdots)$) where $\falpha_k$ (resp.\ $\fbeta_k$) is
a map which to a history $h_k = (x_0, a_0, b_0,
\dots, x_{k-1}, a_{k-1}, b_{k-1}, x_k)$ with $x_i \in \X$, $a_i \in
\A(x_i)$, $b_i\in \B(x_i,a_i)$ for $0\le i\le k$ (resp.\ $(h_k, a_k)$) at stage $k$
associates a probability distribution on a probability space over $\A$
(resp.\ $\B$) which support is included in the possible actions space 
$\A(x_k)$ (resp.\ $\B(x_k, a_k$)). A Markovian (or feedback) strategy
is a strategy which only depends on the information of the current
stage $k$: $\falpha_k$ (resp.\ $\fbeta_k$) depends only on $x_k$
(resp.\ $(x_k, a_k$)), then $\falpha_k(h_k)$ (resp.\ $\fbeta_k(h_k,
a_k)$) will be denoted $\falpha_k(x_k)$ (resp.\ $\fbeta_k
(x_k,a_k)$). It is said stationary if it is independent of $k$, then
$\falpha_k$ is also denoted by $\falpha$ and $\fbeta_k$ by $\fbeta$. A
strategy of any type is said pure if for any stage $k$, the values of
$\falpha_k$ (resp.\ $\fbeta_k$) are Dirac probability measures at
certain actions in $\A(x_k)$  (resp.\ $\B(x_k, a_k$)) then we 
denote also by $\falpha_k$ (resp.\ $\fbeta_k$) the map which to the
history assigns the only possible action in $\A(x_k)$ (resp.\ $\B(x_k,
a_k$)).
 
In particular, if $\salpha$ is a pure Markovian stationary strategy,
then $\salpha = \proc{\falpha_k}$ with $\falpha_k = \falpha$ for all
$k$ and $\falpha$ is a map $\X \rightarrow \A$ such that $\falpha (x)
\in \A(x)$ for all $x\in\X$. In this case, we also speak about pure
Markovian stationary strategy for $\falpha$ and we denote by $\Am$ the
set of such maps. We adopt a similar convention for player \MIN~: $\Bm
:= \set{\bbeta : \X\times\A \rightarrow \B \,|\, \bbeta(x,a) \in \B(x,a)
  \, \forall x \in \X, \, a\in\A(x)}$.

A strategy $\salpha = \proc{\falpha_k}$ (resp.\ $\sbeta =
\proc{\fbeta_k}$) together with 
an initial state determines stochastic processes $\proc{\Ak_k}$ for
the actions of \MAX, $\proc{\Bk_k}$ for the actions of \MIN\ and
$\proc{\Xk_k}$ for the states of the game such that 
\begin{subequations}\label{probG1}
\begin{align}
P(\Xk_{k+1} = y \, | \,  \Ik_k = h_k, \Ak_k = a, \Bk_k = b) & =  p(y\, | \, x,a,b) \\ 
P(\Ak_{k} \in A \, | \,  \Ik_k = h_k) &= \falpha_k(h_k)(A) \\ 
P(\Bk_{k} \in B \, | \,  \Ik_k = h_k, \Ak_k = a) &= \fbeta_k(h_k,a)(B) 
\end{align}
\end{subequations}
where $\Ik_k := (\Xk_0, \Ak_0, \Bk_0, \dots, \Xk_{k-1},
\Ak_{k-1}, \Bk_{k-1} \Xk_k)$ is the history process,
$h_k$ is a history vector at time $k$:
$h_k=(x_0, a_0, b_0, 
\dots, x_{k-1}, a_{k-1}, b_{k-1}, x)$ and  $A$
(resp.\ $B$) are measurable sets in $\A(x)$ ($\B(x, a$) resp.). 
For instance, for each pair of pure Markovian stationary strategies
($\salpha$, $\sbeta$) of the two players, that is such that for $k\ge0$~:
$\falpha_k = \balpha$ with $\balpha \in \Am$ and $\fbeta_k = \bbeta$
with $\bbeta \in \Bm$, the state process  
$\proc{\Xk_k}$ is a Markov chain on $\X$ with transition probability  
\[
P(\Xk_{k+1}=y \, | \,  \Xk_k = x)\,=\, p(y | x, \balpha(x), \bbeta(x, \balpha(x))) \quad \text{ for } x, y \in \X
\]
and  $\Ak_k = \balpha(\Xk_k)$ and $\Bk_k = \bbeta(\Xk_k,\Ak_k)$.

The payoff of the game $\Gamma(x_0)$ starting from $x_0 \in \X$ is the
expected sum of the rewards at all steps of the game 
that \MAX\ wants to maximize and \MIN\ to
minimize. In this paper we consider discounted games $\Gamma_{\mu}$
with discount factor $0 < \mu < 1$: the reward at time $k$ is the payment
made by \MIN\ to \MAX\ times $\mu^k$. When the strategies $\salpha$ for
\MAX\ and $\sbeta$ for \MIN\ are fixed, the payoff of the game
$\Gamma_\mu(x_0, \salpha, \sbeta)$ starting from $x_0$ is then
\[
J(x_0, \salpha, \sbeta) \,=\, \sE^{\salpha, \sbeta}_{x_0} \left[ \, \sum_{k = 0}^{\infty} \mu^k r(\Xk_k,\Ak_k, \Bk_k) \,\right], 
\]
where $\sE^{\salpha, \sbeta}_{x_0}$ denotes the expectation for the
probability law determined by \eqref{probG1}. 
A discounted game can be seen equivalently as a game which has, in
each stage, a stopping probability equal to $1-\mu$, independent of
the actions taken by both players. The value of the game starting from
$x_0 \in \X$, $\Gamma_{\mu}(x_0)$, is then given by 
\begin{equation} \label{valueGdisc}
v(x_0) \,=\, \sup_{\salpha} \inf_{\sbeta} J(x_0, \salpha, \sbeta),
\end{equation}
where the supremum is taken over all strategies $\salpha$ for \MAX\
and the infimum is taken over all strategies $\sbeta$ for \MIN. Note
that a non terminating game without any discount factor (or $\mu = 1$)
is called ergodic. 

We are concerned in finding optimal strategies for both players and
the value of the discounted game $\Gamma_{\mu}$ in each point. These
are given by the dynamic programming equation \cite{Shapley} defined
below. 

\begin{theorem}[Dynamic programming equations \cite{Shapley}]
Assume $\A(x)$  and $\B(x,a)$ are finite sets for all $x\in\X$, $a \in
\A(x)$. Then, the value $v$ of the stochastic game $\Gamma_{\mu}$,
defined in \eqref{valueGdisc}, is the unique solution $v : \X
\rightarrow \R$ of the following dynamic programming equation: 
\begin{equation} \label{eq1}
v(x) \, = \, \underbrace{\max_{a \in \A(x)} \, \left(  \min_{b \in
    \B(x, a)} \, \left( \sum_{y \in \X} \mu \, p(y\, | \, x, a, b) \,
    v(y) \, + \, r(x,a, b) \right)\right)}_{\substack = \, F(v; x)}
    \qquad \forall x \in \X. 
\end{equation}

Moreover, optimal strategies are obtained for both players by taking in~\eqref{valueGdisc}
pure Markovian stationary strategies $\salpha$ for \MAX\ and $\sbeta$ for \MIN\
such that for all $x$ in $\X$, $\balpha(x)$ attains the maximum in
\eqref{eq1}~:
\begin{equation*}
\balpha(x) \,\in\, \argmax{a \in \A(x)} \, F(v; x, a)
\end{equation*}
where
\begin{equation} \label{ILPI}
F(v; x, a) \, := \, \min_{b \in \B(x, a)} \ \left( \underbrace{ \sum_{y \in \X}
\mu \ p(y|x, a, b) \ v(y) \ + \ r(x,a, b)
}_{\substack = \, F(v; x,a,b)}
\right),
\end{equation}
and  for all $x$ in $\X$ and $a$ in $\A(x)$, $\bbeta(x,a)$ attains
the minimum in~\eqref{ILPI}~:
\begin{equation*}
\bbeta(x) \,\in\, \argmin{b \in \B(x,a)} \, F(v; x, a,b) \enspace.
\end{equation*}
Here we use the notation $\operatorname{argmax}_{c \in C} f(c) := \set{c\in C \,\mid\,   f(c) =\max_{c'\in C} f(c')}$ and similarly for $\operatorname{argmin}$.

\end{theorem}

We denote by $F$ the dynamic programming operator from $\R^{\X}$ to itself which maps $v$ to the function 
\begin{equation} \label{DPOperator}
\begin{array}{l l l l}
F(v): & \X & \rightarrow & \R \\
& x & \mapsto & F(v;x)
\end{array}
\end{equation}
where $F(v;x)$ is defined in \eqref{eq1}. This operator is monotone
and contracting with constant $\mu$ in the sup-norm,
i.e.\ $\snorm{F(v) - F(v^{'})}_{\infty} \, \le \, \mu \snorm{v -
  v^{'}}_{\infty}$ for all $v, v^{'} \in \R^\X$.  
Hence, fixed point iterations on Equation~\eqref{eq1}, called value
iterations in the optimal control and game literature, are contracting
for the sup-norm with constant $\mu$.

\section{Two player zero-sum stochastic differential games: the continuous case}
\label{section-SDG}

Another class of games which we consider is the class of two
player differential stochastic games in continuous time. 
In these games, the state space is a regular open
subset $\X$ of $\R^d$ and the dynamics of the game is governed by a
stochastic differential equation which is jointly
controlled by two players (see~\cite{FlemingSoug,swiech} and below). 
In this case, the value of the game
(defined below) is solution of a non linear elliptic partial
differential equation of type~\eqref{eq2intro}, called Isaacs equation (see
also~\cite{FlemingSoug,swiech}). 
The discretization of this equation with a 
monotone scheme in the sense of~\cite{barlsoug} yields the dynamic
programming equation~\eqref{eq1} of a stochastic
game with discrete state space which was described in the previous section.

In the first following subsection, we give the definitions of differential
stochastic games with a bounded state space and a discounted
payoff. Then, in the next subsection, we present a subclass of
these differential games called optimal stopping time games. 
Finally, in the last subsection, we introduce the finite difference 
discretization scheme that we use to discretize the Isaacs
equation~\eqref{BIcont1} and~\eqref{VI} respectively. 
Numerical examples of such kind of games will be presented in
section~\ref{section-numerics}.


\subsection{Differential games with regular controls.}

Assume now that the state space is a regular open subset $\X$ of
$\R^d$. Suppose a probability space $\Omega$ is given, as well as 
a filtration $({\cal F}_t)_{t\geq 0}$ over it (that is a non
decreasing sequence of $\sigma$-algebras over $\Omega$). We consider games
which dynamics is governed by the following stochastic differential
equation~:
\begin{equation} 
d\Xk_t  \,=\,  \bi(\Xk_t,\Ak_t, \Bk_t)\,dt + \sigma (\Xk_t, \Ak_t, \Bk_t)\,dW_{t}, \label{dyn-ar} 
\end{equation}
with initial state $\Xk_0 = x \in \X$. Here $W_{t}$ is a
$d'$-dimensional  Wiener process on $(\Omega, ({\cal F}_t)_{t\geq 0})$;
$\Ak_t$ and $\Bk_t$ are stochastic processes taking values in
closed subsets $\A$ and $\B$ of $\R^p$ and $\R^q$ respectively;
$(x,a,b) \in \X \times \A \times \B \mapsto \bi(x,a,b)
\in \R^d$ and $\X \times \A \times \B \mapsto \sigma (x,
a, b) \in \R^{d\times d'}$ are given functions.
The dimension $d'$ of the Wiener process may be different from $d$ and
is given by the modeling of the problem. 
Assuming that $\Ak_t$ and $\Bk_t$ are adapted to the filtration $({\cal
  F}_t)_{t\geq 0}$ (that is for all $k\ge 0$, $\Ak_t$ and $\Bk_t$ are
${\cal F}_t$-measurable), allows one to define the stochastic process $\Xk_t$
satisfying Equation~\eqref{dyn-ar} and it is a necessary condition to the
assumption that the actions of the two players
depend only on the past states and actions.
We also consider strategies $\salpha =
(\falpha_t)_{t\ge 0}$ (resp.\ $\sbeta = (\fbeta_t)_{t\ge 0}$) of player
\MAX\ (resp.\ \MIN)  
determining the process $(\Ak_t)_{t\ge0}$ (resp.\ $(\Bk_t)_{t\ge0}$). 
In particular, for pure Markovian stationary strategies, one has
$\Ak_t = \balpha (\Xk_t)$ and $\Ak_t = \bbeta (\Xk_t, \Ak_t)$. 

When $\X = \R^d$,  the discounted payoff of the game with discount
rate $ \lambda > 0$ is given by~:
\begin{equation}\label{Jdiff}
 J (x; \salpha, \sbeta) \,=\, \sE^{\salpha, \sbeta}_{x} \left[\,
 \int^{\infty}_{0} e^{-\lambda t} r(\Xk_t,\Ak_t, \Bk_t)\,dt \, | \,
 \Xk_0 =x \,\right] 
\end{equation}
where $(x,a,b) \in \X \times \A \times \B \mapsto r(x,a,b)
\in \R$ is the (instantaneous, or running) reward function.
Now, we consider that $\X$ is a regular open subset $\X $ of $\R^d$.
 In this case, we denote by $\tau$ the first exit time of 
the process $(\Xk_t)_{t\ge0}$ from $\X$, i.e.\ $\tau = \inf \set{t \ge 0 |
  \Xk_t \notin \X}$. Then, the discounted payoff  
of the game stopped at the boundary is~:
\begin{equation}\label{costdiff}
 J (x; \salpha, \sbeta) \,=\, \sE^{\salpha, \sbeta}_{x} \left[\, \int^{\tau}_{0} e^{-\lambda t} r(\Xk_t,\Ak_t, \Bk_t)\,dt + 
e^{- \lambda \tau} \psi_1(\Xk_{\tau}) \, | \, \Xk_0 =x \,\right]
\end{equation}
where the function $ x \in \partial\X \rightarrow \psi_1(x) \in \R$ is called the
terminal reward.
The value function of the differential stochastic game starting from $x$ is defined as in section~\ref{discrete} by 
\begin{equation}\label{valueDiffGame}
v(x) \, = \, \sup_{\salpha} \ \inf_{\sbeta} \ J ( x; \salpha, \sbeta) 
\end{equation}
where the supremum is taken over all strategies $\salpha$ for \MAX\
and the infimum is taken over all strategies $\sbeta$ for \MIN.  

As previously, we are interested in finding the value function of the
game and the corresponding optimal strategies. 
We denote by $L(v;x,a,b)$ the following second order partial differential operator~:
\[ L(v; x, a, b)\, :=\, \sum_{i,j=1}^d \ai_{ij}(x, a, b) \hess{v(x)}{x_i}{x_j}
 +  \sum_{j=1}^d \bi_{j}(x, a,b) \der{v(x)}{x_{j}}  - \lambda v(x) , \]
with ${\displaystyle (\ai_{ij})_{i,j=1,..,d} = \frac{1}{2} \sigma
  \sigma^{T}}$. When $d' \ge d$ and $\sigma(x, a, b)$ is onto for all
$x\in\X,\; a\in \A,\; b\in \B$, 
the matrix $\ai(x, a, b)$ is of
full rank and the operator $L$ is elliptic.
The value of the game $v$ is solution, under some regularity assumptions on
$\Omega$ and on the functions $\bi$, $\sigma$, $r$ and $\psi$
(for instance boundedness and uniform Lipschitz continuity), of the
dynamic programming equation, called Isaacs partial differential
equation~:   
\begin{equation}\label{BIcont1} 
\left\{ \begin{array}{l}
{\displaystyle \max_{a \in \A} \, \left( \min_{b \in \B} \, \left(
 L(v; x,a,b) + r(x,a,b) \,\right) \right) \, =\, 0 \qquad \text{ for
 } x \in \X}
 \\[1em]
v(x)  \,=\,   \psi_1(x)  \qquad \text{ for } x \in  \partial \X.
\end{array} \right. 
\end{equation}
This has been shown in the viscosity sense in \cite{FlemingSoug}. See
also \cite{UserG} and references therein for uniqueness of the
solution of \eqref{BIcont1}. 
If the value $v$ of the game is a classical solution of
\eqref{BIcont1}, $\balpha$ and $\bbeta$ are strategies such that for
all $x$ in $\X$ and $a$ in $\A(x)$, $\balpha(x)$ and $\bbeta(x,a)$ are
the unique actions that realize the maximum and the minimum in
Equation~\eqref{BIcont1} for \MAX\ and \MIN\ respectively, then
$\balpha$ and $\bbeta$ are pure Markovian stationary strategies,
that are optimal for~\eqref{valueDiffGame} (with $\Xk,\Ak,\Bk$ 
satisfying~\eqref{dyn-ar},~\eqref{costdiff}, with
$\Ak_t = \balpha(\Xk_t)$ and $\Bk_t = \bbeta(\Xk_t,\Ak_t)$).

Note that for a game with one player, i.e.\ for a stochastic control problem, Equation \eqref{BIcont1} is the so-called Hamilton-Jacobi-Bellman equation.
Also when $\X$ is bounded, and $L$ is strongly uniformely elliptic (if for some $c>0$, $\ai(x,a,b)\geq c I$ for all $x\in\X,\; a\in \A,\; b\in \B$), then the case $\lambda=0$ can also be considered.

\subsection{Differential games with optimal stopping control}

When the action $(\Ak_t,\Bk_t)$ of the players are not continuous or not bounded, the
dynamic programming equation of the game is no more of the form of
Equation \eqref{BIcont1}, but may be a variational inequality or a
quasi-variational inequality, see for instance~\cite{FriedmanA,Ben78}
for the case of optimal stopping games with one or two
players and ~\cite{FlemingSonerBook,Ben82} for impulse or singular control. 

We consider here an optimal stopping game, that is a game in which one
of the players have the choice of stopping the game at any moment
(see~\cite{FriedmanA} for a more general case). We
assume here that \MAX\ has this ability. Then at each time $t$, he
chooses to stop or not the game, that is he is choosing an element of the action space 
$\{0,1\}$ where $1$ means that the game is continuing, $0$ that the
game stops, with $\Ak_s = 0$ and $\Xk_s = \Xk_t$ for $s \ge t$ when
$\Ak_t = 0$ (i.e.\ $\bi(x, 0, b) = 0$, $\sigma (x, 0, b) = 0$ $\forall
b\in\B, x\in\X$ in~\eqref{dyn-ar}). The second player \MIN\ plays as previously
and we consider the same model as in previous subsection.  
The value of a strategy $\salpha$ for \MAX\ determines a process
$(\Ak_t)_{t\geq 0}$ adapted to the filtration of $(\Xk_t)_{t\ge0}$
(that is $(\sigma(\Xk_t))_{t\ge0}$), then a stopping time
$\stops = \inf\set{t\ge 0 | \Ak_t = 0}$ adapted to the process
$(\Xk_t)_{t\ge0}$ and vice versa.

So if $r(x,0,b) = \lambda \psi_2(x) \,\forall b\in\B$,
 the discounted payoff~\eqref{costdiff} can be written as a function of the
stopping time $\stops$ instead of $\salpha$~:
\[ 
 J (x; \stops, \bbeta) \,=\, \sE^{\stops, \bbeta}_{x} \left[\,
 \int^{\stops}_{0} e^{-\lambda t} r(\Xk_t, 1, \Bk_t)\,dt + e^{- \lambda
 \stops} \psi_2(\Xk_{\stops}) \, \indicator_{\stops < \tau} + e^{- \lambda
 \tau} \psi_1(\Xk_{\tau}) \, \indicator_{\stops = \tau} \, \Big| \, \Xk_0 =x \,\right] \enspace.
\] 
Indeed, if $\stops < \tau$,
then $\Xk_s = \Xk_\stops \in \X$, $s \ge \stops$, so $\tau =
+\infty$, and $ \int_{\stops}^{\tau} e^{-\lambda t} r(\Xk_t, \Ak_t, \Bk_t)\,dt
=  e^{- \lambda  \stops}  \psi_2(\Xk_{\stops})$.
The value function~\eqref{valueDiffGame} of the game starting from $x$
is then given by~:
\[v(x) \, = \, \sup_{\stops} \ \inf_{\bbeta} \ J ( x; \stops, \bbeta) \]
where the supremum is taken over all stopping times $\stops \le \tau$
and the infimum is taken over all strategies $\bbeta$ for \MIN.  

Since the variable ``$a$'' appears only when equal to $1$, one can
ommit it in equations, hence Equation~\eqref{BIcont1} becomes~: 
\begin{equation}\label{VI}
 \left\{ \begin{array}{l}
{\displaystyle
\max \ \bigg\{ \, 
\underbrace{ \min_{b \in \B} \ (\, L(v; x,b) + r(x,b)\,) }_{\substack \unrond }
\, , \, 
\underbrace{\lambda ( \psi_2 (x)  - v (x)) }_{\substack \deuxrond }
\, \bigg\} \,=\, 0 \qquad \text{for $x$ in $\X,$}
}
 \\[1em]
v (x) \,=\,   \psi_1 (x)  \qquad \text{for $x \in \partial\X,$}  
\end{array} \right. 
\end{equation}
since $\lambda > 0$, one can divide the term~\nbcircle{2} by
$\lambda$, and get the variational inequality in the usual form
used in viscosity solutions literature.
In another usual way, Equation~\eqref{VI} can be written
as~:
\begin{equation}
\label{visympli}
\text{for $x \in \X$} \quad
  \left\{
  \begin{array}{l}
  \displaystyle  \min_{b \in \B} \ (\, L(v; x,b) + r(x,b)\,) \le 0  \\ 
    \psi_2(x)  - v(x) \le 0 \\ 
    \left(\displaystyle \min_{b \in \B} \ (\, L(v; x,b) +
    r(x,b)\,)\right) \, \left(\psi_2(x) - v(x)\right)\,=\, 0 
  \end{array}
  \right.
\end{equation}
with $ v (x)=\psi_1 (x)$ for $x \in \partial\X$. Both
Equation~\eqref{VI} and Equation~\eqref{visympli} are called
variational inequalities.
Note however, that Equation~\eqref{VI}, or the resulting equation 
obtained by simplifying by $\lambda$ in~\nbcircle{2}, reveals more the 
control nature and can be used to define viscosity solutions (where one need
to write equations in the form $F(x,v(x), Dv(x), D^2v(x))=0$ on $\X$),
whereas Equation~\eqref{visympli} is more adapted to a variational
approach.

As for \eqref{BIcont1}, if $v$ is a classical solution of~\eqref{VI}
or~\eqref{visympli},
if for all $x$ in $\X$: $\balpha(x)$ is equal to $1$ or $0$ if
resp.\ \nbcircle{1} or \nbcircle{2} is maximum in \eqref{VI}   
and if for all $x$ in $\X$: $\bbeta(x,1)$ is the action $b \in \B$
which realize the minimum in \nbcircle{1},  
then an optimal  pure Markovian stationary strategy is obtained by
taking $\Bk_t = \bbeta(\Xk_t,1)$ and $\stops$ equal to the first
time when $\balpha (\Xk_t) = 0$. 
So this equation behaves as Equation~\eqref{BIcont1} but where the
first player has a discrete action space equal to $\{0,1\}$, $1$
meaning continue to play and $0$ meaning stop the game. This
variational inequality can be treated with the same methods as~\eqref{BIcont1}. 




\subsection{Discretization}
\label{section-discretization}

Several discretization methods may transform equations \eqref{BIcont1}
or \eqref{VI} into a dynamic programming equation of the form
\eqref{eq1}. This is the case when using Markov discrezation of the
diffusion's \eqref{BIcont1} as in \cite{kus77,KuDu92} and in general
when using discrezation schemes that are monotone in the sense of
\cite{barlsoug}. One can obtain such discretizations by using the
simple finite difference scheme below when there are no mixed
derivative (that is $\sigma \sigma^T$ is a diagonal matrix). Under
less restrictive assumptions on the coefficients, finite difference
schemes with larger stencil also lead to monotone schemes
\cite{BoZi03,MuZi05}.  In the deterministic case (when $\sigma \equiv
0$), one can also use semi-Lagrangian scheme \cite{BaFaSo94,BaFaSo99}
or max-plus finite element method \cite{akgbla08}, both of them having
the property of leading to a discrete equation of the
form~\eqref{eq1}.  

We suppose that $\X$ is the $d$-dimensional open unit cube. Let
$h=\frac{1}{m}$ ($m \in{\N}^*$) denote the finite difference step
in each coordinate direction, $e_i$ the unit vector in the
$i^{th}$-coordinate direction, and $x= (x_1, \ldots, x_d)$ a point of
the uniform grid 
$\X_{h}=\X\cap\reseau$. Equation \eqref{BIcont1} is discretized by
replacing the first and second order derivatives of $v$ by the
following approximation, for $i = 1, \dots, d$~: 

\begin{equation}\label{deriv} 
{\partial v(x) \over \partial x_i}
\sim {v(x + he_i) - v(x - he_i) \over 2h}
\end{equation} 
or
\begin{equation}\label{decent} 
{\partial v(x) \over \partial x_i} \sim
\left\{ \begin{array}{lll} 
\displaystyle {v(x + he_i) - v(x) \over h} & \text{when} & \bi_i(x,a,b) \geq 0 \\ \\ 
\displaystyle {v(x) - v(x-he_i) \over h} & \text{when} & \bi_i(x,a,b) < 0. 
\end{array} \right. 
\end{equation}
\begin{equation}  
\hes{v}{x_i}(x)  \sim   { {v(x+ h e_i )  -2v(x)+ v(x- h
e_i )} \over {h}^2 }, \label{centre}
\end{equation} 
Approximation \eqref{deriv} may be used when $L$ is uniformly elliptic and $h$ is small, whereas \eqref{decent} has to be used when $L$ is degenerate  (see \cite{kus77,KuDu92}).
For equations \eqref{BIcont1} and \eqref{VI}, these differences are computed in the entire grid $\X_h$, by prolonging $v$ on the ``boundary'', $\partial \X_h := \partial \X \cap$\reseau
using Dirichlet boundary condition: 
\[
v(x) \,=\, \psi(x) \qquad \forall\, x \in \pX\cap\reseau.
\]

We obtain a system of $N_h$ non linear equations of $N_h$ unknowns,
the values of the function $v_h: x \in \X_h \mapsto v_h(x) \in \R$~:  
\begin{equation}\label{eq5.4}
\max_{a \in \A} \, (\, \min_{b \in \B} \ (\, L_h(v_h; (x,a,b)) + r(x,a,b) \,)\,) \,=\, 0 \;\;\;
\forall\, x\in\X_h \enspace,
\end{equation} 
where  $N_h=\sharp \X_h\sim 1/h^d$ and $L_h$ is a function which to $v
\in \R^{\X_h}$, $x \in \X_h$, $a \in \A$, $b \in \B$ associates the
approximation of $L(v;x,a,b)$. 

When there are no mixed derivatives ($\ai_{i,j} (x,a,b) = 0$ if $i\ne
j$, $i,j \in \{1, \dots, d\}$), the discretization is monotone in the
sense of \cite{barlsoug}, then if \eqref{BIcont1} has a unique
viscosity solution, the solution $v_h$ of~\eqref{eq5.4} converges
uniformly to the solution $v$ of~\eqref{BIcont1}~\cite{barlsoug}.  
Moreover, multiplying Equation~\eqref{eq5.4} by $ch^2$ with $c$ 
small enough, it  can be rewritten in the form~\eqref{eq1},
with a discount factor $\mu=1 - O(\lambda c h^2)$.
A similar result holds for the discretization of~\eqref{VI}
(by multiplying only the diffusion part by $ch^2$).

We refer to section~\ref{section-expl_Ieq} for an example of an Issacs
equation~\eqref{Iex1} whose discretization (using
scheme~\eqref{decent}-~\eqref{centre}) yields an
equation~\eqref{Iex1d} which has the form of~\eqref{eq1}.

\section{Background for numerical solution of discrete dynamic programming equations}
\label{section-backG}

In this section, we present the policy iteration algorithm to
solve the dynamic programming equation~\eqref{eq1} of a two player
zero-sum discounted stochastic game with finite state space. 
We first present the policy iteration algorithm for a one  player game
which is then used in the following subsection to define the policy
iteration algorithm for the two player case. 
The last part of this section is devoted to a recall of multigrid
methods which we will use in the policy iterations for solving the
linear systems.

\subsection{Policy iteration algorithm for one player games}

\label{sectionPI1P}

First, we consider a one player stochastic game with a \MIN\
player and finite state space $\X$. In this case, the dynamic
programming operator $F$, mapping 
$\R^n$ to itself, is given for each $x\in\X$ by~:
\begin{equation} \label{eq1opdef}
F(v; x) \ = \ \min_{b \in \B(x)} \ \left( \sum_{y \in \X}
\mu \ p(y|x, b) \ v(y) \ + \ r(x, b)\right) \enspace.
\end{equation}
This game is more commonly called a Markov Decision Process (MDP) with
finite state space $\X$, we refer to~\cite{Howard60,Denardo-Fox68,puterman} for a deeper description on this topic.
Then, the discounted value of the game starting in $x \in\X$ is given
by~:
\[
v(x) \,=\, \inf_{\sbeta} \sE^{\sbeta}_{x} \left[ \, \sum_{k =
    0}^{\infty} \mu^k r(\Xk_k,\Bk_k) \,\right], 
\]
where the processes $\Xk_k,\Bk_k$ and strategies $\sbeta$ are defined
such as in the section~\ref{discrete}. The value $v$ of the game is
solution of the dynamic  
programming equation~: $v(x) \, = \, F(v;x)$ for $x$ in $\X$.
Then the policy iteration algorithm for Markov Decision Processes, that was
first introduce by Howard~\cite{Howard60}, is given in
Algorithm~\ref{algo-MDP} and give us the discounted value of the game
$v: \X \rightarrow \R$ and the optimal policy for \MIN. 

\begin{algorithm}\caption{Policy iteration algorithm for Markov
    Decision Processes (one player game)}
\label{algo-MDP}
\noindent Given an initial policy $\bbeta^0\in\Bm$, the
policy iterations consist in applying successively the two following
steps:
\begin{enumerate}
\item \label{step-val-MDP} 
Compute the value $v^{\kj+1}$ of the game with fixed feedback
policy $\bbeta^\kj$, that is the solution  of 
\begin{equation} \label{MDP-LS}
v^{\kj+1} (x) = \sum_{y \in \X} \mu \ p(y|x, \bbeta^\kj(x)) \ v(y) \ + \ r(x, \bbeta^\kj(x))
\end{equation}
\item \label{step-improve-MDP} 
Improve the policy: Find the optimal feedback policy
$\bbeta^{\kj+1}$ for the value $v^{\kj+1}$, i.e.\ for each $x$ in $\X$,
chose $\bbeta^{\kj+1}(x)$ such that~:
$$
 \bbeta^{\kj+1}(x) \ \in \ \underset{b \in
  \B(x)}{\operatorname{argmin}} \ \left\{\, \sum_{y \in \X}
\mu \, p(y\, | \, x,  b) \, v^{\kj+1}(y) \, + \, r(x, b)\, \right\} 
$$

\end{enumerate}
until we cannot improve the policy anymore.
\end{algorithm}

Each policy iteration of Algorithm~\ref{algo-MDP} strictly improves the
current policy and produces a non increasing sequence of values
$(v^{\kj})_{\kj\geq 1}$. It implies that the algorithm never visits twice
the same policy. Hence if the action sets are finite in each point of
$\X$, the policy iterations stop after a finite time (see for 
instance~\cite{PutBrum,LionsMercier80,Bertsekas87}).
Moreover, under regularity assumptions, the policy 
iteration algorithm for a one player game with infinite action spaces
is equivalent to Newton's
method~\cite{aki90b,BankRose82,zidaniBokanowski09,PutBrum}.
Indeed, define $G(v) = F(v) - v$, then the problem is to find the
solution of $G(v) = 0$ where all entries of $G$ are concave
functions.  The policy improvement step can be seen as the computation
of an element of the sup-differential of $G$ in the current
approximation $v^{\kj+1}$ and the value improvement step computes the zero
of the previous sup-differential. When $G$ is regular, the sequence of
value functions $(v^\kj)_{\kj\ge 1}$ is exactly the sequence of the
Newton's algorithm.

\subsection{Policy iteration algorithm for two player games}
\label{section-PI2P}

Now, we  give the policy iteration algorithm for solving a two player
zero-sum stochastic games with finite state space $\X$, as defined in
Puri thesis~\cite{Puri}. 
Recall the definitions of section~\ref{discrete}, we need to solve the
dynamic programming Equation~\eqref{eq1} which give us the value of
the game (Equation~\eqref{valueGdisc}) and the optimal
strategies for both players.
For a fixed pure feedback policy for \MAX\ $\balpha\in \Am$, the value $v$ of the game is solution of the
equation $v = F^{\balpha}(v)$ where $F^{\balpha}$ is an operator mapping
$\R^n$ to itself whose $x$-coordinate is given by~:
\[ 
F^{\balpha} (v;x) \, := \, F(v;x,\balpha(x)) \, = \, \min_{b \in \B(x, \balpha(x))} \ \left( \sum_{y \in \X}
\mu \ p(y|x, \balpha(x), b) \ v(y) \ + \ r(x,\balpha(x), b) \right) \enspace,
\] 
for each $x \in\X$ and $v \in \R^n$.
Note that $F^{\balpha}$ is the dynamic programming
operator of a one player game with only the \MIN\ player.
Then the policy iteration algorithm is given in
Algorithm~\ref{PIalgo}.

\begin{algorithm}
\caption{Policy Iteration}
\label{PIalgo}
\noindent Given an initial policy $\balpha^0\in \Am$ for \MAX, the
policy iterations consist in applying successively the two following
steps:
\begin{enumerate}
\item \label{step-value-PI} Compute the value $v^{\ki+1}$ of the game with fixed feedback
policy $\balpha^\ki$, that is the solution  of 
\[
v^{\ki+1} \ = \ F^{\balpha^\ki}(v^{\ki+1})
\]
by using Algorithm~\ref{algo-MDP}.
\item Improve the policy: Find the optimal feedback policy
$\balpha^{\ki+1}$ of \MAX\ for the value $v^{\ki+1}$ , i.e.\ for each $x$ in $\X$,
chose $\balpha^{\ki+1}(x)$ such that~:
\[
\balpha^{\ki+1}(x) \ \in \ \underset{a \in
  \A(x)}{\operatorname{argmax}} \ F (v^{\ki+1};x,a) 
\]
where $F (v;x,a)$ is defined by~\eqref{ILPI}.
\end{enumerate}
until we cannot improve the policy anymore.
\end{algorithm}

Step~\ref{step-value-PI} of Algorithm~\ref{PIalgo} is performed by
using the policy iteration algorithm for a one player game. That is,
given an initial feedback policy for \MIN\ $\bbeta^{\ki,0} \in \Bm$, we
iterate on  
\MIN\ policies $\bbeta^{\ki,\kj}\in \Bm$ and value functions $v^{\ki,\kj}$.  
Then at each step $\kj$ of the interior policy iteration
(Algorithm~\ref{algo-MDP} step~\ref{step-val-MDP}), one computes 
$v^{\ki,\kj+1}$, the value of the game with fixed strategies
$\balpha^\ki \in \Am$
for \MAX\ and $\bbeta^{\ki,\kj} \in \Bm$ for \MIN. 
This is done by solving the linear system~: 
\begin{equation} \label{LSPI}
v^{\ki,\kj+1} \ = \ \mu \, M^{\balpha^{\ki} \bbeta^{\ki,\kj}} \,
v^{\ki,\kj+1} \,+\, r^{\balpha^{\ki} \bbeta^{\ki,\kj}}\enspace,
\end{equation}
where for all $\balpha\in \Am, \,\bbeta\in \Bm$: $M^{\balpha\bbeta}
\in \R^{n \times n}$ is a 
stochastic matrix whose elements are defined by
$(M^{\balpha\bbeta})_{x,y} = p(y|x, \balpha(x),
\bbeta(x))$ for all $x,y \in \X$ and $r^{\balpha\bbeta} \in \R^n$ is
the vector whose elements are defined 
by $(r^{\balpha \bbeta})_x = r(x, \balpha(x)
\bbeta(x))$ for $x \in\X$.


As for the one player case, each iteration of the policy iteration
algorithm strictly improve the current policy, hence it
can never visit twice the same policy. Moreover, the algorithm
produces a non decreasing (resp.\ non increasing) sequence of values
$(v^\ki)_{\ki\geq 1}$  (resp.\ $(v^{\ki,\kj})_{\kj\geq 1}$) of the
external loop (resp.\ internal loop), see~\cite{Puri,CochetGaub}.
It follows that if the action sets for both players are finite in each point of
$\X$, the policy iterations stop after a finite time~\cite{Puri}.

\subsection{AMG}

The linear systems defined in~\eqref{LSPI} have all the form $v =
\mu M v + r$ where $M$ a Markov matrix. We solve them using
algebraic multigrid methods which we recall in this section. 

Standard multigrid was originally created in the seventies to solve
efficiently linear elliptic partial differential equations (see for
instance~\cite{MR972752}). It works as follows. 
Multigrid methods require discretizations of the given continuous
equation on a sequence of grids.
Each of them, starting from a coarse grid, being a refinement of the
previous until a given accuracy is attained. The size of the coarsest
grid is chosen such that the cost of solving the problem on it is
cheap. Assume also that transfer operators between these grids are
given: interpolation and restriction. Then, a multigrid cycle on the
finest grid consists in~: first, the application of a smoother on the
finest grid; then a restriction of the residual on the next coarse
grid; then solving the residual problem on this coarse grid using the
same multigrid scheme; then, interpolate this solution (which is an
approximation of the error) and correct the error on the fine grid; 
finally, the application of a smoother on the finest grid.
If the multigrid components are properly chosen, this process is
efficient to find the solution on the finest grid. Indeed, in general
the relaxation process is smoothing the error which then can be well
approximated by elements in the range of the interpolation.  
It implies, in good cases, that the contraction factor of the
multigrid method is independent of the discretization step and also
the complexity is in the order of the number of discretization
points. We shall refer to this standard method as geometric multigrid.

Algebraic multigrid method, called AMG, has been initially developed
in the early eighties (see for example~\cite{Brandt,Brandt1,stub1})
for solving large sparse linear systems arising from the
discretization of partial differential equations with unstructured
grids or PDE's not suitable for the application of the geometric
multigrid solver or large discrete problems not derived from any
continuous problem.  

The AMG method consists of two phases, called ``setup phase'' and
``solving phase''. In contrast to geometric multigrids, the mode of
constructing the coarse levels (coarse ``grids'') which constitute the
setup phase, is based only on the algebraic equations. The points of
the fine grids are represented by the variables and coarse grids by
subset of these variables.  The selection of those coarse variables
and the construction of the transfer operators between levels are done
in such a way that the range of the interpolation approximates the
errors not reduced by a given relaxation scheme.  Then the ``solving
phase'' is performed in the same way as a geometric multigrid method
and consists of the application of a smoother and a correction of the
error by a coarse grid solution. The whole process is briefly recall
below. 

Consider a system of $n$ linear equations given in the matrix form:
\begin{equation}\label{LS}
  Av = f  
\end{equation}
where the matrix $A$ $\in \R^{n\times n}$ and the vector $f$ $\R^n$ are given, and we are looking for the vector $v \in \R^n$. We call fine grid $\Omega^0$ the set of all variables of the system, i.e.\ $\Omega^0 = \{1, \ldots, n\}$.

First, recall that a relaxation method consist of the following approximations:
\[ u \leftarrow S u + \So f   \quad  \text{with}  \quad S = I- \So A\]
where $S$ is called the smoothing operator and $I$ is the identity operator in $\R^{n\times n}$. The error $e = u - v$ propagates as
\[ e \leftarrow S e. \]
The method is said to converge if $\rho(S) < 1$ where $\rho(S) = \max_i |\lambda_i|$ is the spectral radius of $S$ with $\lambda_i$ his eigenvalues. 
For example, the smoother operator of the weighted Jacobi method is $S = I-wD^{-1}A$ and that of the Gauss-Seidel is $S = I-L^{-1}A$ where $D$ and $L$ are the diagonal and lower triangular part of the matrix $A$ resp.\

Assume $\Omega^l$ the grid on level $l$ where level $0$ correspond to
the finest grid $\Omega^0$. The construction of the coarse grid
$\Omega^{l+1}$ from the fine grid $\Omega^l$, consists in the
splitting of the $n_l$ variables from the grid $\Omega^l$ into two
distinct subsets, namely $C$ which contains the variables belonging to
both grids, $\Omega^l$ and $\Omega^{l+1}$, and $F$ the variables
belonging to the grid $\Omega^l$ only. We have then $\Omega^l = C \cup
F$. The coarse grid $\Omega^{l+1} = C$ contains $n_{l+1}$
variables. This splitting is based on the ``connections'' between the
variables on level $l$~\cite{Brandt,stub1} and such as the range of
the associate interpolation or prolongation operator $\Pl^l_{l+1}$
accurately approximates the errors not efficiently reduced by the
relaxation phase (these errors are ``smooth'' in the algebraic
multigrid terminology). The restriction operator $\Rl^{l+1}_l$ maps
residuals from grid $\Omega^l$ to the grid $\Omega^{l+1}$. 
In~\cite{Brandt,stub1}, the operator is fixed to be $\Rl^{l+1}_l =
(\Pl_{l+1}^l)^T$. The coarse grid operator is defined by $A^{l+1} =
\Rl^{l+1}_lA^l\Pl^l_{l+1}$ where $A^{l+1}$ is the approximation of
$A^l$ on $\Omega^{l+1}$ and $A^0 = A$. Similarly, for any vector $v^l
\in \R^{n_l}$ we denote $v^{l+1} = \Rl^{l+1}_l v^l$ its restriction on
$\Omega^{l+1}$. This construction can be repeated recursively from the
finest level $l=0$ to the coarsest level $L$. 

The solution phase consists in applying the multigrid cycle described
in Algorithm~\ref{MGcycle}, it is called V($\nu_1$,$\nu_2$)-cycle if
$\gamma = 1$ and W($\nu_1$,$\nu_2$)-cycle if $\gamma = 2$. 
\begin{algorithm}
\caption{Multigrid scheme $u^l \leftarrow MG(u^l, f^l)$}
\label{MGcycle}
\begin{algorithmic}
\IF{ $l < L$} 
\STATE {\bf pre relaxation~:} 
\STATE $\qquad u^l \leftarrow S u^l + \So f^l \quad$ (on $\Omega^l$) $\qquad \nu_1$ times 
\STATE {\bf coarse grid correction~:} 
\STATE $\qquad f^{l+1} \leftarrow \Rl^{l+1}_l (f^l - A^l u^l)$
\STATE $\qquad u^{l+1} \leftarrow 0$
\STATE $\qquad u^{l+1} \leftarrow MG(u^{l+1}, f^{l+1}) \qquad \gamma$ times
\STATE $\qquad u^l \leftarrow u^l + \Pl^l_{l+1} u^{l+1}$
\STATE {\bf post relaxation~:} 
\STATE $\qquad u^l \leftarrow S u^l + \So f^l \quad$ (on $\Omega^l$) $\qquad \nu_2$ times 
\ELSE  
\STATE Solve $A^L u^L = f^L$
\ENDIF
\end{algorithmic}
\end{algorithm}
Convergence theorem for the V-cycle is given in\cite{stub1} for
$A$ symmetric and positive definite. See also~\cite{Brandt,Brandt1,Falgout}, for two-level convergence for linear systems
where the matrix of the system is a M-matrix, 
symmetric and positive definite. 
Also we can find in the literature, two-grid convergence analysis for
non-symmetric linear system in~\cite{Notay} and~\cite{Mense}.

\section{A multigrid algorithm for discrete dynamic programming equations}
\label{section-AMGpi}

\subsection[Policy iteration combined with algebraic multigrid method]{Policy iteration combined with algebraic multigrid method (AMG$\pi$)}

\begin{figure}
\setlength{\unitlength}{1973sp}%
\begingroup\makeatletter\ifx\SetFigFont\undefined%
\gdef\SetFigFont#1#2#3#4#5{%
  \reset@font\fontsize{#1}{#2pt}%
  \fontfamily{#3}\fontseries{#4}\fontshape{#5}%
  \selectfont}%
\fi\endgroup%
\begin{picture}(10582,4699)(676,0)
\thinlines
{\color[rgb]{0,0,0}\put(2326,4200){\line( 0,-1){1950}}
}%
\put(2476,4000){\makebox(0,0)[lb]{\smash{{\SetFigFont{10}{12.0}{\rmdefault}{\mddefault}{\updefault}{\color[rgb]{0,0,0}$\balpha^0$}%
}}}}
\put(2476,3500){\makebox(0,0)[lb]{\smash{{\SetFigFont{10}{12.0}{\rmdefault}{\mddefault}{\updefault}{\color[rgb]{0,0,0}$\vdots$}%
}}}}
\put(2476,3100){\makebox(0,0)[lb]{\smash{{\SetFigFont{10}{12.0}{\rmdefault}{\mddefault}{\updefault}{\color[rgb]{0,0,0}$\balpha^\ki$}%
}}}}
\put(2476,2600){\makebox(0,0)[lb]{\smash{{\SetFigFont{10}{12.0}{\rmdefault}{\mddefault}{\updefault}{\color[rgb]{0,0,0}$\vdots$}%
}}}}
{\color[rgb]{0,0,0}\put(5626,3400){\line( 0,-1){1950}}
}%
\put(5776,3100){\makebox(0,0)[lb]{\smash{{\SetFigFont{10}{12.0}{\rmdefault}{\mddefault}{\updefault}{\color[rgb]{0,0,0}$\bbeta^{\ki,0}$}%
}}}}
\put(5776,2600){\makebox(0,0)[lb]{\smash{{\SetFigFont{10}{12.0}{\rmdefault}{\mddefault}{\updefault}{\color[rgb]{0,0,0}$\vdots$}%
}}}}
\put(5776,2200){\makebox(0,0)[lb]{\smash{{\SetFigFont{10}{12.0}{\rmdefault}{\mddefault}{\updefault}{\color[rgb]{0,0,0}$\bbeta^{\ki,\kj}$}%
}}}}
\put(5776,1700){\makebox(0,0)[lb]{\smash{{\SetFigFont{10}{12.0}{\rmdefault}{\mddefault}{\updefault}{\color[rgb]{0,0,0}$\vdots$}%
}}}}
{\color[rgb]{0,0,0}\put(9301,2600){\line( 0,-1){2000}}
}%
\put(9451,2300){\makebox(0,0)[lb]{\smash{{\SetFigFont{10}{12.0}{\rmdefault}{\mddefault}{\updefault}{\color[rgb]{0,0,0}$v^{\ki,\kj,0}$}%
}}}}
\put(9451,1800){\makebox(0,0)[lb]{\smash{{\SetFigFont{10}{12.0}{\rmdefault}{\mddefault}{\updefault}{\color[rgb]{0,0,0}$\vdots$}%
}}}}
\put(9451,1500){\makebox(0,0)[lb]{\smash{{\SetFigFont{10}{12.0}{\rmdefault}{\mddefault}{\updefault}{\color[rgb]{0,0,0}$v^{\ki,\kj,\km}$}%
}}}}
\put(9451,1000){\makebox(0,0)[lb]{\smash{{\SetFigFont{10}{12.0}{\rmdefault}{\mddefault}{\updefault}{\color[rgb]{0,0,0}$\vdots$}%
}}}}
\put(9451,700){\makebox(0,0)[lb]{\smash{{\SetFigFont{10}{12.0}{\rmdefault}{\mddefault}{\updefault}{\color[rgb]{0,0,0}$v^{\ki,\kj+1,0}$}%
}}}}
{\color[rgb]{0,0,0}\put(6876,2200){\vector( 1, 0){2000}}
}%
{\color[rgb]{0,0,0}\put(3201,3100){\vector( 1, 0){2000}}
}%
\put(676,3200){\makebox(0,0)[lb]{\smash{{\SetFigFont{10}{12.0}{\rmdefault}{\mddefault}{\updefault}{\color[rgb]{0,0,0}PI external}%
}}}}
\put(4076,2400){\makebox(0,0)[lb]{\smash{{\SetFigFont{10}{12.0}{\rmdefault}{\mddefault}{\updefault}{\color[rgb]{0,0,0}PI intern}%
}}}}
\put(8226,1600){\makebox(0,0)[lb]{\smash{{\SetFigFont{10}{12.0}{\rmdefault}{\mddefault}{\updefault}{\color[rgb]{0,0,0}AMG}%
}}}}
\end{picture}%
\caption{Representation of the nested iterations of AMG$\pi$.}
\label{amgPIscheme}
\end{figure}
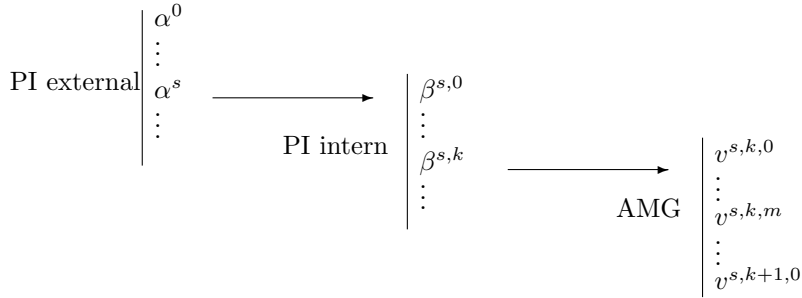

Recall that in the policy iteration algorithm for games at each step
$\kj$ of the interior policy iteration, we have to solve a linear
system~\eqref{LSPI} which is of the form $ v =  \mu M v + r$
with $M$ a Markov matrix and $0<\mu<1$ the discount factor. Since
$(I - \mu M)$ are non singular $M$-matrices, we use AMG to solve those
systems. For shortness in the sequel, we shall  
call the resulting algorithm AMG$\pi$ that is the combination of
policy iterations and AMG. The name AMG$\pi$ refers also to the
numerical implementation of this algorithm. 
Note that in practice, in Algorithm~\ref{algo-MDP} (equivalently in
Algorithm~\ref{PIalgo}), the policy iterations are stopped when after
Step~\ref{step-val-MDP}, the norm of the residual,   
$r_v = F(v) - v$, is smaller than a given value denoted 
by $\epsilon$. We used this stopping criterion in AMG$\pi$.
The iterations of AMG$\pi$  are summarized in the scheme represented in 
Figure \ref{amgPIscheme} 
\noindent where $(v^{\ki,\kj,0}, \cdots, v^{\ki,\kj,\km}, \cdots
,v^{\ki,\kj+1,0})$ is a sequence of value functions generated by the
multigrid solver. 
The algebraic multigrid methods allows us to solve linear systems
arising from either the discretization of Isaacs or
Hamilton-Jacobi-Bellman equations or a true finite state space
zero-sum two player game.
However in the present paper, we restrict ourselves to numerical tests
for the discretization of stochastic differential games, since the AMG
algorithm needs some improvements to be applied to arbitrary non
symmetric linear systems arising in game problems. 
 

In the one player game case, convergence results of combination of
policy iteration and geometric multigrid method have been established
by Hoppe~\cite{Hope86,Hope87} and Akian~\cite{akian,aki90b}.

\subsection[Full multi-level policy iteration]{Full multi-level policy iteration (FAMG$\pi$)}
\label{sectionFullamgPI}

Recall that the number of policy iterations can be exponential in the
cardinality of the state space $\X$. 
However, as for Newton's algorithm, convergence can be improved by starting the
policy iterations with a good initial guess, close to the solution.  
With this in mind, we present a full multi-level scheme, that we shall
call FAMG$\pi$. 
As in standard FMG, starting from the coarsest
level, it consists in solving the problem at each grid
level by performing policy iterations AMG$\pi$ until a convergence
criterion is verified, then to interpolate the strategies and value
function to
the next level, in order to initialize the policy iterations of that
level.   
This scheme is repeated until the finest level is attained. 

\begin{figure}
\center
\begin{picture}(0,0)%
\includegraphics{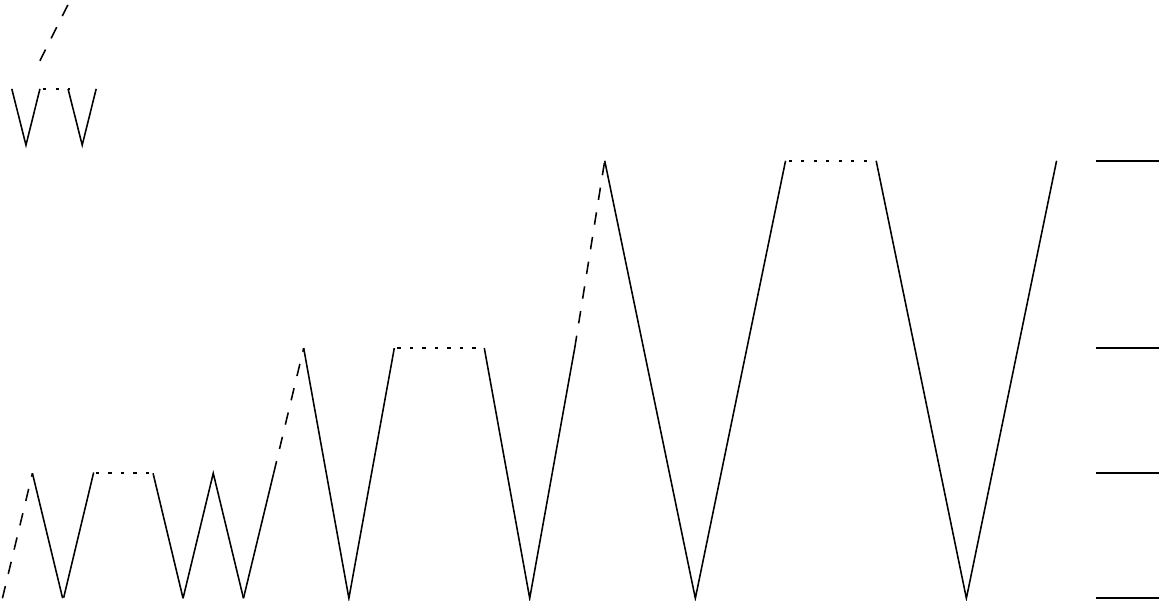}
\end{picture}%
\setlength{\unitlength}{3947sp}%
\begingroup\makeatletter\ifx\SetFigFont\undefined%
\gdef\SetFigFont#1#2#3#4#5{%
  \reset@font\fontsize{#1}{#2pt}%
  \fontfamily{#3}\fontseries{#4}\fontshape{#5}%
  \selectfont}%
\fi\endgroup%
\begin{picture}(7523,3024)(664,-3748)
\put(1301,-1050){\makebox(0,0)[lb]{\smash{{\SetFigFont{10}{12.0}{\rmdefault}{\mddefault}{\updefault}{\color[rgb]{0,0,0}{Interpolation of strategies and value}}%
}}}}
\put(1301,-1450){\makebox(0,0)[lb]{\smash{{\SetFigFont{10}{12.0}{\rmdefault}{\mddefault}{\updefault}{\color[rgb]{0,0,0}{AMG$\pi$}}%
}}}}
\put(6376,-1661){\makebox(0,0)[lb]{\smash{{\SetFigFont{10}{12.0}{\rmdefault}{\mddefault}{\updefault}{\color[rgb]{0,0,0}$\X^0$}%
}}}}
\put(6376,-2561){\makebox(0,0)[lb]{\smash{{\SetFigFont{10}{12.0}{\rmdefault}{\mddefault}{\updefault}{\color[rgb]{0,0,0}$\X^{1}$}%
}}}}
\put(6376,-3161){\makebox(0,0)[lb]{\smash{{\SetFigFont{10}{12.0}{\rmdefault}{\mddefault}{\updefault}{\color[rgb]{0,0,0}$\X^{2}$}%
}}}}
\put(6376,-3761){\makebox(0,0)[lb]{\smash{{\SetFigFont{10}{12.0}{\rmdefault}{\mddefault}{\updefault}{\color[rgb]{0,0,0}$\X^3$}%
}}}}
\end{picture}%
\caption{FAMG$\pi$ with AMG$\pi$ V-cycles}
\label{FMGamgPI}
\end{figure}
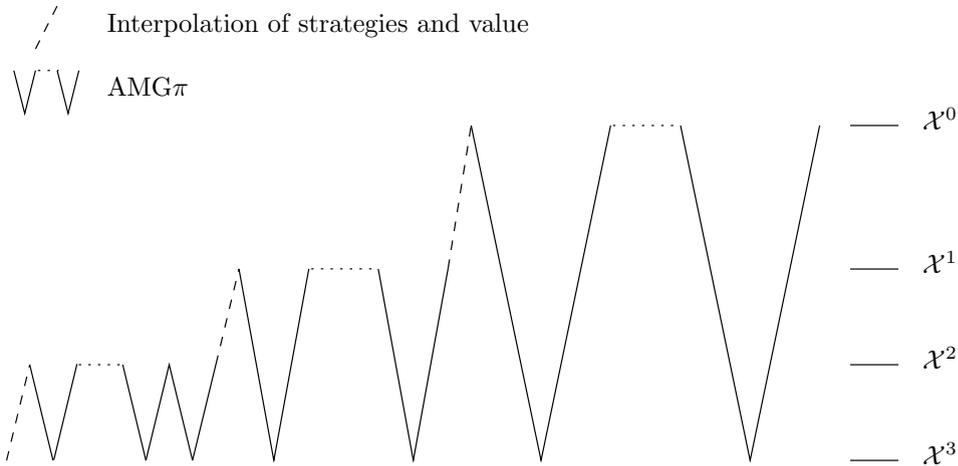

The algorithm FAMG$\pi$ only applies to Isaacs partial differential
equations~\eqref{BIcont1}. It works as follows. The state
space $\X$ is first discretized on sequence of $L_F + 1$ grids~:
$\X_{L_F} \subset \dots \subset \X_{1} \subset \X_0 = \X_h$ such that on grid
$\X_l$, $0\le l\le L_F$, the discretization step is $h_l = 2^l h$, 
where $h$ is the discretization step chosen on the finest grid
$\X_h$. Then, the Isaacs PDE is discretized on all levels, $0\le l\le
L_F$, using the finite differences scheme~\eqref{decent}-~\eqref{centre}. 
For level $l$, we denote by $F_l : \X_l \rightarrow \X_l$ the dynamic
programming operator, $(v)^l : \X_l
\rightarrow \R$ the value of game, $x \in \X_l 
\rightarrow (\balpha)^l(x) \in \A(x)$ and
$(x \in \X_l, a\in \A(x)) \rightarrow (\bbeta)^l(x,a) \in \B(x,a)$ the
strategies of \MAX\ and \MIN\ respectively.
We denote by $\Il^{l-1}_{l}$ the linear
interpolation operator which maps any vector $(v)^l$ from $\R^{\X_{l}}$ to 
$\R^{\X_{l-1}}$~:
\[
\Il^{l-1}_{l}(v)^l(x) = \quad
  \left\{
  \begin{array}{l l}
    (v)^l(x) &  x \in \X_l \\
    \sum_{y \in \Nx(x)} \frac{1}{\sharp(\Nx)} \, (v)^l(y) & x \in
    \X_{l-1}\setminus\X_l 
  \end{array}
  \right.
\]
where $\Nx(x) = \set{y \in \X_l
  \,|\, \norm{x-y}_2 <= h_l}$ for $x \in \X_{l-1}\setminus\X_l$,
and we denote by $\Ul^{l-1}_{l}$ the operator which
interpolates a 
strategy from grid $\X_{l}$ to grid $\X_{l-1}$, for
instance for a strategy of \MAX~: 
\[
\Ul^{l-1}_{l}((\balpha)^l)) = \quad
  \left\{
  \begin{array}{l l}
    (\balpha)^l(x) &  x \in \X_l \\
    a_0 \in  \A(x) & x \in
    \X_{l-1}\setminus\X_l 
  \end{array}
  \right.
\]
where $a_0$ is chosen arbitrary $A(x)$ in for $x \in
\X_{l-1}\setminus\X_l$.
We denote by AMG$\pi(\balpha, \bbeta, v, \epsilon)$ 
the algorithm AMG$\pi$ with initial strategy $\balpha$ for player
\MAX\ iterations, initial policy $\bbeta$ for the first iteration of
player \MIN, value $v$
as initial approximation for the first call of AMG and $\epsilon$ the stopping
criterion for the policy iterations.
Then FAMG$\pi$ algorithm is given in Algorithm~\ref{algo-FAMGPI} where
$c >  0$ is a given constant.

\begin{algorithm} 
\caption{FAMG$\pi$} 
\label{algo-FAMGPI}
\begin{algorithmic}
\STATE Given an initial $(\balpha^0)^{L_F}, (\bbeta^0)^{L_F}$ and $(v^0)^{L_f}$ on level $L_F$,
\FOR{$l = L_F$ to $1$}
\STATE  
$((\balpha)^l,(\bbeta)^l, (v)^l) \leftarrow$ AMG$\pi((\balpha^0)^l,
(\bbeta^0)^l, (v^0)^l, c h_l^2)$ on level $l$
\STATE $(v^0)^{l-1} \, = \, \Il^{l-1}_{l} \, (v)^l$
\STATE $(\alpha^0)^{l-1} = \Ul^{l-1}_{l}(\alpha)^{l}$ and $(\beta^0)^{l-1} = \Ul^{l-1}_{l}(\beta)^{l}$
\ENDFOR
\STATE solve $v = F (v)$ on $\X_h$ by using
AMG$\pi((\balpha^0)^0, (\bbeta^0)^0, (v^0)^0, \epsilon)$
\end{algorithmic}

\end{algorithm}

Figure \ref{FMGamgPI} illustrates the FAMG$\pi$
algorithm when V-cycles are use in AMG$\pi$. The dashed lines
represent the interpolation of the solution and strategies from a
coarse grid $\X^{l}$ to the next fine grid $\X^{l-1}$. The
continuous V-lines are the V-cycles of AMG$\pi$ which are not fixed in
number since at each level, AMG$\pi$ cycles are performed until a given
criterion is attained. 

Note that our FAMG$\pi$ program only applies to stochastic
differential games since for them coarse representation, including
equations and strategies, can be easily constructed by tacking
different sizes of discretization step. 

For one-player discounted games with infinite number of actions and
under regularity and strong convexity assumptions, it is shown
in~\cite{aki90b,akian} 
that this kind of full multi-level policy iteration has a computing time 
in the order of the cardinality of $\X$.

\section{Numerical results}
\label{section-numerics}

In this section, we apply our programs AMG$\pi$ and FAMG$\pi$, 
which were implemented in C, to examples of two player zero-sum
stochastic differential games. 
Let first give some details about the implementation of the algorithms
that we use and some notations for the numerical results.

The AMG linear solver of AMG$\pi$ implements the construction phase,
including the 
coarsing scheme and the interpolation operator, described
in~\cite{stub1} and the general recursive multigrid cycle for the 
solution phase (see Algorithm~\ref{MGcycle}). In the tests,
W($1$,$1$)-cycles were used and the chosen smoother is a CF relaxation
method, that is a Gauss Seidel relaxation scheme that relaxes first on
C-points and then on F-points.   
The AMG$\pi$ program is the implementation of the method explained in
section~\ref{section-AMGpi} with the above AMG linear solver. 
The FAMG$\pi$ program is the implementation of
Algorithm~\ref{algo-FAMGPI}.   

The following notations are used in the tables: $\ki$ denotes the
iteration over \MAX\ policies and $\nkj$ is 
the corresponding number of iterations for \MIN\ policies, that is the
number of linear systems solved at iteration $\ki$.  
The residual error of the game is denoted by $r_v = F(v) - v$ and the exact
error, when known, by $e = F(v) - u$ where $u$ is the discretized
exact solution of the game. 
The infinite norm and discrete $L_2$ norm are given for each of them. 


\subsection{Isaacs equations}
\label{section-expl_Ieq}

\begin{figure} 
\center
\includegraphics[height=5cm]{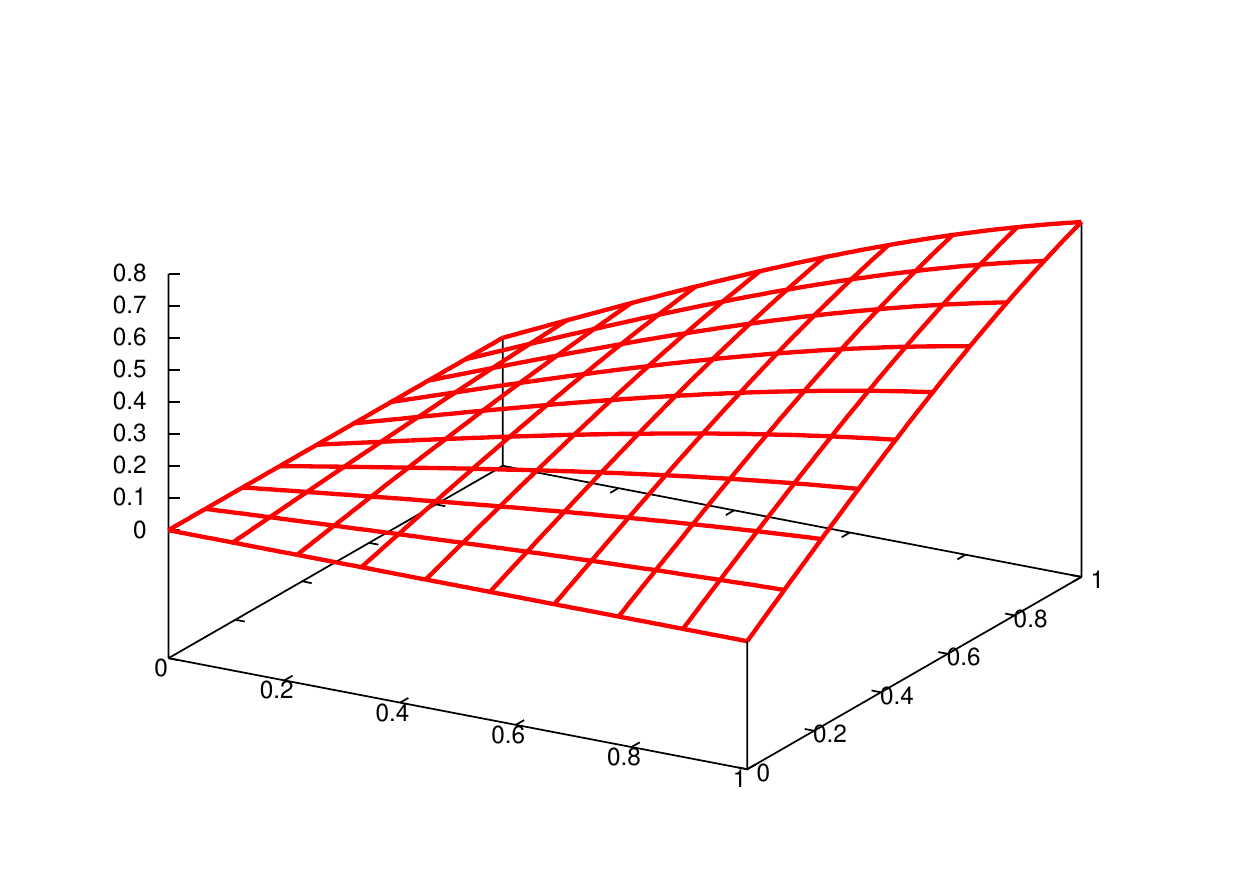}
\caption{Graph of $\sin(x_1) \times \sin(x_2)$ on $\X = [0,1] \times [0,1]$.}
\label{sinxsiny}
\end{figure}

The first example concern a diffusion problem where the value
$v:\X\rightarrow\R$ of 
the game is solution of the following Isaacs PDE~: 
\begin{equation}\label{Iex1}
\left\{
\begin{array}{l l}
{\displaystyle 
\max_{a \in \A} \min_{b \in \B} \left(
\Delta v(x) + (a \cdot \nabla v (x)) - \left(b \cdot \nabla v (x)
\right) - \lambda v(x) + \frac{\norm{b}^2_2}{2} + f (x)
 \right) = 0 }
& x \text{ in } \X \enspace, \\
 v(x) = \psi_1(x) & x \text{ in } \partial\X 
\end{array}
\right.
\end{equation}
where $\X = ]0,1[ \times ]0,1[$ is the unit square, $\A =
\set{a \in \R^2 \,|\, \norm{a}_2 \le 1}$, $\B = \R^2$,
$\psi_1(x_1,x_2) = \sin(x_1) \times \sin(x_2)$ for $(x_1,x_2) \in \partial\X$, and $f(x)
= - (\Delta u(x) + \norm{\nabla u(x)}_2 - 0.5 \norm{\nabla
u(x)}_2^2 - \lambda u(x))$ with
$u(x_1,x_2) = \sin(x_1) \times \sin(x_2)$ for $x = (x_1,x_2) \in \X$. 
Note that the exact solution is $v(x_1,x_2) = \sin(x_1) \times
\sin(x_2)$ on $\X = [0,1] \times [0,1]$ and is represented in
Figure~\ref{sinxsiny}. Indeed, by convex duality (or computation of
Fenchel-Legendre transformations~\cite{ROCK}), we
have that 
\[\norm{u}_2 = \max_{{\norm{a}}_2\le 1, a \in \R^d}\, a \,\cdot\, \,u  \qquad
    \text{ and } \qquad \frac{1}{2}\norm{u}^2_2 =
    \max_{b \in \R^d}\, b\,\cdot\,u - \frac{1}{2} \norm{b}^2_2\]
for all $u \in \R^d$, $a = \frac{u}{\norm{u}_2^2}$ and $b=u$ are
optimal solutions in these equations. 

To solve Equation~\eqref{Iex1}, we first discretize the domain
$[0,1]\times[0,1]$ on a grid with $m+1$ 
points in each direction, i.e.\ with a discretization step 
$h=\frac{1}{m}$ and we obtain a discrete space $\X_h$ with boundary $
\partial \X_h$. We denote by $x_i = i h$ with $i=0,\dots,m$ such that
$\X_h = \set{(x_i,x_j) \,|\, i,j \in \set{1,\dots,m-1}}$ and $\partial
\X_h = \set{(x_i,x_j) \,|\,   i \in  \set{0,m}, j
\in \set{0,\dots,m} \text{ or } j \in \set{0,m}, i \in \set{0,\dots,m}}$.
Then, using the discretization  
scheme~\eqref{decent}-~\eqref{centre}, Equation~\eqref{Iex1} becomes 
for $(x_i,x_j) \in \X_h$~:
\begin{align*}
0\, = & \max_{(a_1,a_2) \in \A} \min_{(b_1,b_2) \in \B} \ \Bigg\{
 \left( \frac{- 4 v(x_i,x_j)  + v(x_{i+1},x_j) + v(x_{i-1},x_j) +
  v(x_i,x_{j+1}) + v(x_i,x_{j-1})}{h^2}  \right) \\
  & +  (a_1 - b_1) \left( \frac{v(x_{i+1},x_j) -  v(x_i,x_j)}{h}
\right) \indicator_{(a_1 - b_1)\ge 0}
+ (a_1 - b_1) \left( \frac{v(x_i,x_j) - v(x_{i-1},x_j)}{h}
\right) \indicator_{(a_1 - b_1) < 0}  \\
 & + (a_2 - b_2) \left( \frac{v(x_{i},x_{j+1}) -  v(x_i,x_j)}{h}
\right) \indicator_{(a_2 - b_2) \ge 0} 
 + (a_2 - b_2) \left( \frac{v(x_i,x_j) - v(x_{i},x_{j-1})}{h}
\right) \indicator_{(a_2 - b_2) < 0}  \\
 & - \lambda v(x_i,x_j) + \frac{b_1^2 + b_2^2}{2} + f(x_i,x_j)
 \Bigg\} \enspace,
\end{align*}
multiply by $\frac{h^2}{c}$, where $c= 4 + h \abs{a_1 - b_1} + h \abs{a_2 - b_2} > 0$, and adding
$v(x_i,x_j)$ on both sides, we obtain~:
\begin{align} 
v(x_i,x_j) \, = & \max_{(a_1,a_2) \in \A} \min_{(b_1,b_2) \in \B} \
  \left(1 + \frac{h^2}{c} \lambda \right) ^{-1} \nonumber\\
 &  \Bigg\{   
 \left( \frac{1}{c} + \frac{h}{c} (a_1 - b_1) \indicator_{(a_1 -  b_1)\ge
 0} \right) v(x_{i+1},x_j) 
 +  \left( \frac{1}{c} - \frac{h}{c} (a_1 - b_1) \indicator_{(a_1 - b_1) < 0} \right) v(x_{i-1},x_j) \nonumber\\
 & +  \left( \frac{1}{c} + \frac{h}{c} (a_2 - b_2) \indicator_{(a_2 -  b_2)\ge
 0} \right) v(x_{i},x_{j+1}) 
 +  \left( \frac{1}{c} - \frac{h}{c} (a_2 - b_2) \indicator_{(a_2 -
 b_2) < 0} \right) v(x_{i},x_{j-1}) \nonumber \\
 & + \frac{h^2}{c} \,\frac{b_1^2 + b_2^2}{2} + \frac{h^2}{c} \,
 f(x_i,x_j)  \ \Bigg\} \qquad \text{ for } (x_i,x_j) \in \X_h
 \enspace, 
\label{Iex1d}
\end{align}
where $v(x_i,x_j)$ is replaced by $ \psi_1(x_i,x_j)$ for
$(x_i,x_j)\in\partial\X_h $. 
This equation has the form of Equation~\eqref{eq1} with a discount
factor $\mu$ equal to $(1 + \frac{h^2}{c} \lambda) 
^{-1} \le 1$, transition probabilities from $(x_i,x_j) \in \X_h$ to
$(x_{i^{'}},x_{j^{'}}) \in \X_h$ are given by~: 
\begin{align}
p ((x_{i^{'}},x_{j^{'}}) | (x_i,x_j),
(a_1,a_2), (b_1,b_2)) \,=\, 
&\frac{1}{c} + \frac{h}{c} (a_1 - b_1) \indicator_{(a_1 - b_1)\ge 0} 
&  \text{ if } i^{'}=i+1,j^{'}= j \enspace, \label{Iex1pr} \\
&\frac{1}{c} - \frac{h}{c} (a_1 - b_1) \indicator_{(a_1 - b_1) < 0} 
& \text{ if } i^{'}=i-1,j^{'}= j \enspace, \nonumber \\
&\frac{1}{c} + \frac{h}{c} (a_2 - b_2) \indicator_{(a_2 - b_2)\ge  0} 
& \text{ if } i^{'}=i,j^{'}= j+1 \enspace,  \nonumber \\
&\frac{1}{c} - \frac{h}{c} (a_2 - b_2) \indicator_{(a_2 - b_2) < 0} 
&  \text{ if } i^{'}=i,j^{'}= j-1 \enspace, \nonumber \\
&0 &  \text { else} \enspace, \nonumber 
\end{align}
and the running cost is, for $(x_i,x_j)\in\X_h$~: 
\begin{align*}
r((x_i,x_j), (a_1,a_2), (b_1,b_2))\, =&\, \frac{h^2}{c} \left(\frac{b_1^2 + 
b_2^2}{2} +  f(x_i,x_j)\right) 
 \\
 & + \left(  \frac{h}{c} (a_1 - b_1) \indicator_{(a_1 -  b_1)\ge
 0} \right) \psi_1(x_{i+1},x_j)  \indicator_{(x_{i+1},x_j) \in
  \partial \X_h} \\
 & -  \left( \frac{h}{c} (a_1 - b_1) \indicator_{(a_1 - b_1) < 0}
\right) \psi_1(x_{i-1},x_j) \indicator_{(x_{i-1},x_j)\in 
 \partial \X_h} \\
 & +  \left( \frac{h}{c} (a_2 - b_2) \indicator_{(a_2 -  b_2)\ge
 0} \right) \psi_1(x_{i},x_{j+1}) \indicator_{(x_{i},x_{j+1}) \in\partial\X_h}\\
 & -  \left( \frac{h}{c} (a_2 - b_2) \indicator_{(a_2 - b_2) < 0}
 \right) \psi_1(x_{i},x_{j-1}) \indicator_{(x_{i},x_{j-1}) 
 \in\partial\X_h}\enspace. 
\end{align*}
Note that when $i,j\in\{2,\dots,m-2\}$ the sum of
the transition probabilities from $(x_i,x_j)$ to the points of $\X_h$
equals $\mu$, when $i$ or $j$ is in $\{1,m-1\}$ this sum
is strictly less than $\mu$. Hence, the matrix $M^{\balpha,\bbeta}$ in
\eqref{LSPI} is substochastic, and since it is irreducible, it has a
spectral radius strictly less than one. So even when $\lambda = 0$ or
equivalently $\mu = 1$, the system~\eqref{LSPI} has an unique solution
and the dynamic programing equations has also an unique solution. 
Hence, we shall take $\lambda = 0$ in the numerical tests. 
Note also that for this example, the matrices $M^{\balpha,\bbeta}$ in
\eqref{LSPI} are 
not symmetric but close to be symmetric when $h$ is small, since the
non-symmetric part correspond to the order one term in
equation~\eqref{Iex1d} and are dominated by order two terms when $b$
is optimal in \eqref{Iex1d}.

\begin{table}
\centering
\caption{Numerical results for equation~\eqref{Iex1} on a $1025\times
  1025$ points grid.}
\label{table1Iex1}
\begin{tabular}{| c| c| c| c| c| c| c| c| c|}
\multicolumn{7}{c}{Policy iteration with LU} \\
\hline 
 $\ki$ & $\nkj$ &  $\norm{r_v}_{\infty}$ &
 $\norm{r_v}_{L_2}$ &  $\norm{e}_{\infty}$ &
$\norm{e}_{L_2}$ & cpu time (s) \\
\hline 
$ 1 $ & $ 3 $ & $ 8.51e-7 $ & $ 5.96e-7 $ & $ 4.47e-2 $ & $ 2.48e-2 $ & $ 1.40e+2 $ \\
\hline 
$ 2 $ & $ 2 $ & $ 2.44e-8 $ & $ 6.16e-9 $ & $ 1.84e-4 $ & $ 1.05e-4 $ & $ 2.31e+2 $ \\
\hline 
$ 3 $ & $ 1 $ & $ 7.38e-13 $ & $ 2.03e-13 $ & $ 4.13e-6 $ & $ 2.16e-6 $ & $ 2.77e+2 $ \\
\hline 
\end{tabular}
\ \\
\begin{tabular}{| c| c| c| c| c| c| c| c| c|}
\multicolumn{7}{c}{AMG$\pi$} \\
\hline 
$\ki$ & $\nkj$ &  $\norm{r_v}_{\infty}$ &
 $\norm{r_v}_{L_2}$ &  $\norm{e}_{\infty}$ &
$\norm{e}_{L_2}$ & cpu time (s) \\
\hline 
$ 1 $ & $ 3 $ & $ 8.51e-7 $ & $ 5.96e-7 $ & $ 4.47e-2 $ & $ 2.48e-2 $ & $ 2.65e+1 $ \\
\hline 
$ 2 $ & $ 2 $ & $ 2.44e-8 $ & $ 6.16e-9 $ & $ 1.84e-4 $ & $ 1.05e-4 $ & $ 4.59e+1 $ \\
\hline 
$ 3 $ & $ 1 $ & $ 7.92e-13 $ & $ 2.02e-13 $ & $ 4.13e-6 $ & $ 2.16e-6 $ & $ 5.56e+1 $ \\
\hline 
\end{tabular}
\end{table}

In tables~\ref{table1Iex1}, we present numerical results when
equations~\eqref{Iex1} is discretized on a grid with $1025$ points in
each direction, i.e.\ with a discretization step of $h = 1/2^{10}$. 
The stopping criterion for the policy iterations is $\epsilon = 10^{-10}$. 
The first table of~\ref{table1Iex1} shows the results of the policy
iteration algorithm with a direct solver LU (we used the package
UMFPACK~\cite{UMFPACK04}) and the second table of~\ref{table1Iex1} the
results of AMG$\pi$.  
We observe that AMG$\pi$ solves the problem faster than
the policy iterations with a direct solver. 
In both tables, we see that only three steps on \MAX\ policies
are needed (first column) and a total of six steps on \MIN\ policies
(second column) which involves the resolution of six linear systems.
The small number of iterations is due to the fact that the solution is
regular.
In table~\ref{table-IexFMG}, we show that the computation time is
improved when applying FAMG$\pi$ with $c=0.1$ to the same example. In
this case, the problem is solved in approximately $18s$.

\begin{table} 
\centering
\caption{Numerical results for Equation \eqref{Iex1} on a $1025
  \times 1025$ points grid, computed by FAMG$\pi$ with $c=10^{-1}$.}
\label{table-IexFMG}
\begin{tabular}{| c| c| c| c| c| c| c|}
\hline 
 $\ki$ & $\nkj$ &  $\norm{r_v}_{\infty}$ &  $\norm{r_v}_{L_2}$ &  $\norm{e}_{\infty}$ & $\norm{e}_{L_2}$ & cpu time (s) \\
\hline 
\multicolumn{7}{|c|}{ points in each direction : $ 3 $, h $ 5.00e-01 $ } \\ 
\hline 
$ 1 $ & $ 2 $ & $ 1.42e-01 $ & $ 1.42e-01 $ & $ 1.07e-01 $ & $ 1.07e-01 $ & $ << 1 $ \\
\hline 
$ 2 $ & $ 1 $ & $ 2.34e-03 $ & $ 2.34e-03 $ & $ 2.45e-04 $ & $ 2.45e-04 $ & $ << 1 $ \\
\hline 
\multicolumn{7}{|c|}{ points in each direction : $ 5 $, h $ 2.50e-01 $ } \\ 
\hline 
$ 1 $ & $ 2 $ & $ 5.53e-03 $ & $ 2.84e-03 $ & $ 3.00e-03 $ & $ 1.75e-03 $ & $ << 1 $ \\
\hline 
\multicolumn{7}{|c|}{ points in each direction : $ 9 $, h $ 1.25e-01 $ } \\ 
\hline 
$ 1 $ & $ 2 $ & $ 2.40e-04 $ & $ 1.10e-04 $ & $ 8.20e-04 $ & $ 4.46e-04 $ & $ << 1 $ \\
\hline 
\multicolumn{7}{|c|}{ points in each direction : $ 17 $, h $ 6.25e-02 $ } \\ 
\hline 
$ 1 $ & $ 2 $ & $ 3.18e-05 $ & $ 7.83e-06 $ & $ 3.36e-04 $ & $ 1.90e-04 $ & $ 1.00e-02 $ \\
\hline 
\multicolumn{7}{|c|}{ points in each direction : $ 33 $, h $ 3.12e-02 $ } \\ 
\hline 
$ 1 $ & $ 1 $ & $ 5.89e-04 $ & $ 7.08e-05 $ & $ 5.05e-04 $ & $ 1.99e-04 $ & $ 1.00e-02 $ \\
\hline 
\multicolumn{7}{|c|}{ points in each direction : $ 65 $, h $ 1.56e-02 $ } \\ 
\hline 
$ 1 $ & $ 1 $ & $ 1.69e-04 $ & $ 1.25e-05 $ & $ 1.62e-04 $ & $ 4.67e-05 $ & $ 4.00e-02 $ \\
\hline 
\multicolumn{7}{|c|}{ points in each direction : $ 129 $, h $ 7.81e-03 $ } \\ 
\hline 
$ 1 $ & $ 1 $ & $ 4.28e-05 $ & $ 2.16e-06 $ & $ 4.73e-05 $ & $ 1.21e-05 $ & $ 1.80e-01 $ \\
\hline 
\multicolumn{7}{|c|}{ points in each direction : $ 257 $, h $ 3.91e-03 $ } \\ 
\hline 
$ 1 $ & $ 1 $ & $ 1.08e-05 $ & $ 3.77e-07 $ & $ 1.31e-05 $ & $ 6.07e-06 $ & $ 7.50e-01 $ \\
\hline 
\multicolumn{7}{|c|}{ points in each direction : $ 513 $, h $ 1.95e-03 $ } \\ 
\hline 
$ 1 $ & $ 1 $ & $ 2.70e-06 $ & $ 6.61e-08 $ & $ 7.29e-06 $ & $ 3.56e-06 $ & $ 3.13e+00 $ \\
\hline 
\multicolumn{7}{|c|}{ points in each direction : $ 1025 $, h $ 9.77e-04 $ } \\ 
\hline 
$ 1 $ & $ 2 $ & $ 1.23e-10 $ & $ 8.13e-13 $ & $ 4.16e-06 $ & $ 2.17e-06 $ & $ 1.85e+01 $ \\
\hline 
\end{tabular}
\end{table}

\begin{figure} 
\input{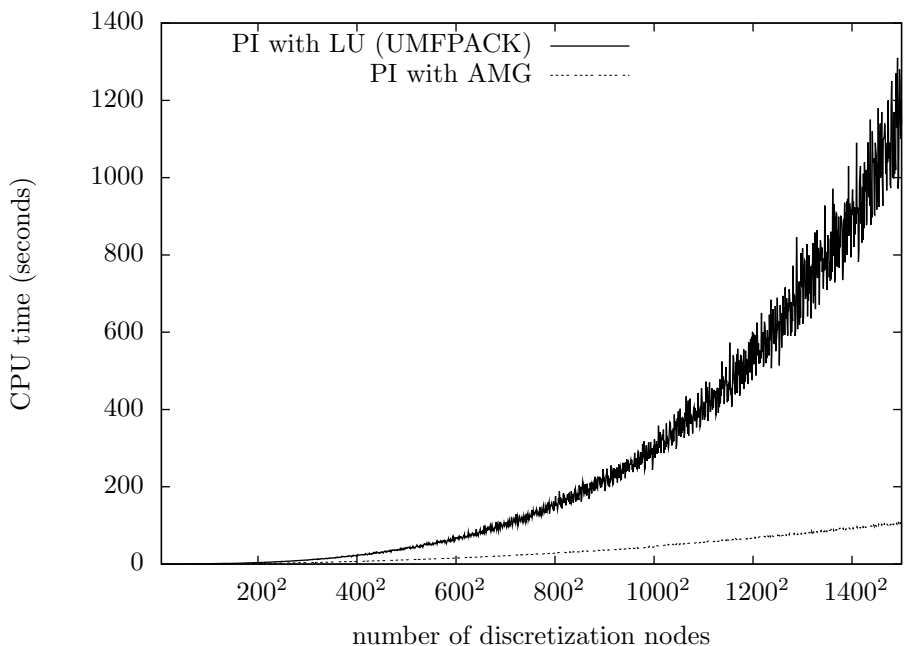}
\caption{Comparison between AMG$\pi$ versus policy iteration algorithm
  with a LU
  solver for solving equation~\eqref{Iex1} when increasing the size of
  the problem.}
\label{timeVSsize}
\end{figure}

\begin{figure} 
\input{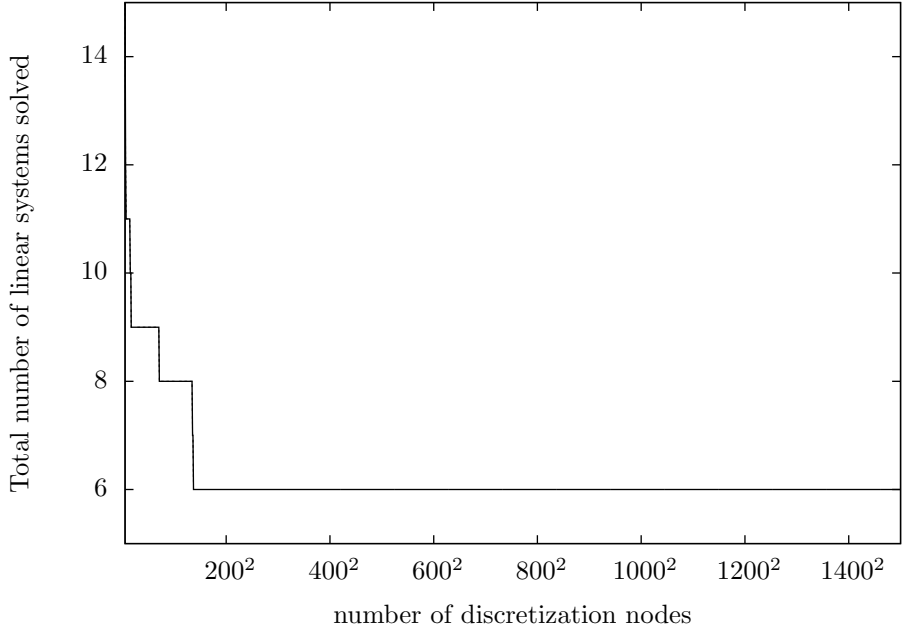}
\caption{Number of iterations on \MIN\ policies (i.e the number of
  linear systems solved) for solving equation~\eqref{Iex1} when
  increasing the size of the problem corresponding to
  figure~\ref{timeVSsize} for both methods (AMG$\pi$ and policy iteration algorithm with LU).}
\label{timeVSsizeIt}
\end{figure}

In figure~\ref{timeVSsize}, we compare the policy iteration
algorithm with a direct solver LU (UMFPACK~\cite{UMFPACK04}) and
AMG$\pi$ for solving equation~\eqref{Iex1d}, when
increasing by one the 
number of discretization points in each direction from $m = 5$ to $m =
1500$. 
The stopping criterion for the policy iterations is $\epsilon = 10^{-10}$.
In figure~\ref{timeVSsizeIt}, we represent the 
corresponding number of iterations on \MIN\ policies, i.e the number
of linear systems solved for each size of problem, this number is
the same for both methods. 
We can see that the most part of the computation time for the
resolution of the non-linear equation~\eqref{Iex1d} is used to solved
the linear systems involved in the policy iteration. 
We also remark that the computation time for AMG$\pi$ seems to grow
linearly with the size of the problem.

\begin{table} 
\centering
\begin{tabular}{| c | c| c| c| c| c| c| c|} 
\hline 
 $\ki$ & $\nkj$ & AMG &  $\norm{r_v}_{\infty}$ &  $\norm{r_v}_{L_2}$ &  $\norm{e}_{\infty}$ & $\norm{e}_{L_2}$ & cpu time (s) \\ 
\hline 
$1$  & $2$  &  $ 5 , 4  $    & $2.15e-04$  & $1.52e-04$  & $4.45e-02$  & $2.50e-02$  & $5.00e-02$  \\ 
 \hline 
$2$  & $2$  &  $4 , 3  $    & $5.97e-06$  & $1.59e-06$  & $2.36e-04$  & $1.43e-04$  & $1.00e-01$  \\ 
 \hline 
$3$  & $1$  &  $ 3   $    & $3.02e-09$  & $7.47e-10$  & $6.49e-05$  & $3.44e-05$  & $1.30e-01$  \\ 
 \hline 
\end{tabular} 

\caption{Numerical results with a $65 \times 65$ points grid, computed by AMG$\pi$ for equation \eqref{Iex1}.}\label{Iex1t1} 
\begin{tabular}{| c | c| c| c| c| c| c| c|} 
\hline 
 $\ki$ & $\nkj$ & AMG &  $\norm{r_v}_{\infty}$ &  $\norm{r_v}_{L_2}$ &  $\norm{e}_{\infty}$ & $\norm{e}_{L_2}$ & cpu time (s) \\ 
\hline 
$1$  & $2$  &  $ 5 , 4  $    & $5.40e-05$  & $3.80e-05$  & $4.46e-02$  & $2.49e-02$  & $2.30e-01$  \\ 
 \hline 
$2$  & $2$  &  $4 , 3  $    & $1.53e-06$  & $3.95e-07$  & $2.07e-04$  & $1.23e-04$  & $4.30e-01$  \\ 
 \hline 
$3$  & $1$  &  $3  $    & $4.08e-10$  & $9.65e-11$  & $3.28e-05$  & $1.72e-05$  & $5.40e-01$  \\ 
 \hline 
\end{tabular} 

\caption{Numerical results with a $129 \times 129$ points grid, computed by AMG$\pi$ for equation \eqref{Iex1}.}\label{Iex1t2} 
\begin{tabular}{| c | c| c| c| c| c| c| c|} 
\hline 
 $\ki$ & $\nkj$ & AMG &  $\norm{r_v}_{\infty}$ &  $\norm{r_v}_{L_2}$ &  $\norm{e}_{\infty}$ & $\norm{e}_{L_2}$ & cpu time (s) \\ 
\hline 
$1$  & $2$  &  $5 , 4  $    & $1.35e-05$  & $9.51e-06$  & $4.47e-02$  & $2.49e-02$  & $1.06e+00$  \\ 
 \hline 
$2$  & $2$  &  $4 , 3  $    & $3.86e-07$  & $9.86e-08$  & $1.94e-04$  & $1.13e-04$  & $1.98e+00$  \\ 
 \hline 
$3$  & $1$  &  $ 3  $    & $5.17e-11$  & $1.22e-11$  & $1.65e-05$  & $8.63e-06$  & $2.49e+00$  \\ 
 \hline 
\end{tabular} 

\caption{Numerical results with a $257 \times 257$ points grid, computed by AMG$\pi$ for equation \eqref{Iex1}.}\label{Iex1t3} 
\begin{tabular}{| c | c| c| c| c| c| c| c|} 
\hline 
 $\ki$ & $\nkj$ & AMG &  $\norm{r_v}_{\infty}$ &  $\norm{r_v}_{L_2}$ &  $\norm{e}_{\infty}$ & $\norm{e}_{L_2}$ & cpu time (s) \\ 
\hline 
$1$  & $2$  &  $  5 , 4 $    & $3.39e-06$  & $2.38e-06$  & $4.47e-02$  & $2.48e-02$  & $4.55e+00$  \\ 
 \hline 
$2$  & $2$  &  $ 4 , 3   $    & $9.71e-08$  & $2.46e-08$  & $1.87e-04$  & $1.08e-04$  & $8.28e+00$  \\ 
 \hline 
$3$  & $1$  &  $ 3 $    & $6.26e-12$  & $1.55e-12$  & $8.26e-06$  & $4.31e-06$  & $1.04e+01$  \\ 
 \hline 
\end{tabular} 

\caption{Numerical results with a $513 \times 513$ points grid, computed by AMG$\pi$ for equation \eqref{Iex1}.}\label{Iex1t4} 
\begin{tabular}{| c | c| c| c| c| c| c| c|} 
\hline 
 $\ki$ & $\nkj$ & AMG &  $\norm{r_v}_{\infty}$ &  $\norm{r_v}_{L_2}$ &  $\norm{e}_{\infty}$ & $\norm{e}_{L_2}$ & cpu time (s) \\ 
\hline 
$1$  & $2$  &  $ 5 , 4 $    & $8.48e-07$  & $5.95e-07$  & $4.47e-02$  & $2.48e-02$  & $1.85e+01$  \\ 
 \hline 
$2$  & $2$  &  $  4 , 3  $    & $2.43e-08$  & $6.15e-09$  & $1.83e-04$  & $1.05e-04$  & $3.40e+01$  \\ 
 \hline 
$3$  & $1$  &  $ 3 $    & $7.40e-13$  & $2.02e-13$  & $4.13e-06$  & $2.16e-06$  & $4.27e+01$  \\ 
 \hline 
\end{tabular} 

\caption{Numerical results with a $1025 \times 1025$ points grid, computed by AMG$\pi$ for equation \eqref{Iex1}.}\label{Iex1t5} 
\begin{tabular}{| c | c| c| c| c| c| c| c|} 
\hline 
 $\ki$ & $\nkj$ & AMG &  $\norm{r_v}_{\infty}$ &  $\norm{r_v}_{L_2}$ &  $\norm{e}_{\infty}$ & $\norm{e}_{L_2}$ & cpu time (s) \\ 
\hline 
$1$  & $2$  &  $ 5 , 4 $    & $2.12e-07$  & $1.49e-07$  & $4.47e-02$  & $2.48e-02$  & $7.46e+01$  \\ 
 \hline 
$2$  & $2$  &  $  4 , 3  $    & $6.09e-09$  & $1.54e-09$  & $1.82e-04$  & $1.04e-04$  & $1.38e+02$  \\ 
 \hline 
$3$  & $1$  &  $ 3 $    & $1.13e-13$  & $3.04e-14$  & $2.07e-06$  & $1.08e-06$  & $1.72e+02$  \\ 
 \hline 
\end{tabular} 

\caption{Numerical results with a $2049 \times 2049$ points grid, computed by AMG$\pi$ for equation \eqref{Iex1}.}\label{Iex1t6}

\end{table}

Each table \ref{Iex1t1} to \ref{Iex1t6} contains numerical results
for Equation~\eqref{Iex1} discretized  on grids with discretization
step $h = \frac{1}{2^6}$, $h = \frac{1}{2^7}$,  $h = \frac{1}{2^8}$, $h =
\frac{1}{2^9}$, $h = \frac{1}{2^{10}}$ and $h = \frac{1}{2^{11}}$ respectively. 
For these tests, the stopping criterion for the policy iterations is
$\epsilon = 0.001\,{h}^2$ where $h$ is the discretization step. The
stopping criterion for the linear solver AMG is $\norm{r}_2 < 10^{-12}$
where $r$ is the residual for the linear system.
For each line of the tables, the third column, named AMG, contains the number of
iterations needed by AMG for solving each linear system ($\nkj$
systems per line). 
We can see that the number of iterations of AMG is independent
of the size of the problem. 
Note that the norm of the error $\norm{e}$ decrease slowly when the
grid becomes finer, this is because the exact solution
(Figure~\ref{sinxsiny}) is smooth and a small number of points is
sufficient to get a good approximation, also the non-linearity of the
problem gives a worse approximation than one might expect in the linear case.
But a smooth solution is generally more difficult for linear
iterative solvers.

\subsection{Optimal stopping game}

\begin{figure}
\center
\includegraphics[height=5cm]{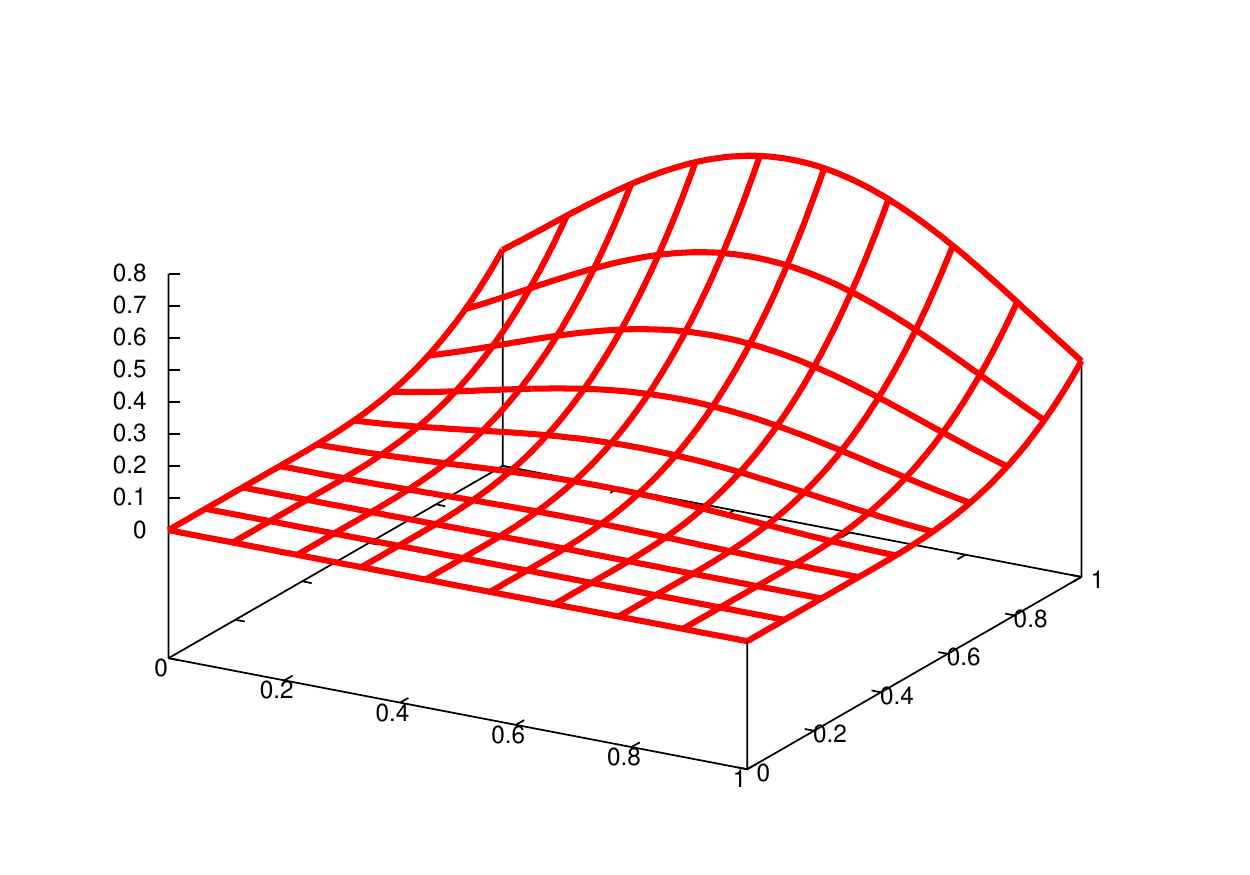}
\caption{Graph of the solution of equation \eqref{VIex1}.}
\label{solvipGraphSol}
\end{figure}

Next tests concern an optimal stopping time game where the value
$v:\X\rightarrow\R$ of the game is solution of the variational
inequality~: 
\begin{equation}\label{VIex1}
\left\{
\begin{array}{l l}
{\displaystyle 
\max_{a \in \A} \left\{
 \underbrace{
\min_{b \in \B} \left(
0.5 \Delta v(x) - \left(b \cdot \nabla v (x)
\right) + \frac{\norm{b}^2_2}{2} + f (x)
 \right)
}_{\substack \unrond}
 \,,\, \underbrace{\psi_2(x) -
    v(x)}_{\substack \deuxrond} \right\} \,=\, 0 }
& x \text{ in } \X \\
v(x) = \psi_1(x) & x \text{ in } \partial\X 
\end{array}
\right.
\end{equation}
where $\X = ]0,1[ \times ]0,1[$ is the unit square, the sets 
$\A = \set{0, 1}$, $\B = \R^2$, $\psi_2(x_1,x_2)=0$ for
$(x_1,x_2) \in \X$, for $(x_1, x_2) \in \X$~:
\[
f (x_1, x_2) \, = \,\left\{
\begin{array}{l l}
-(0.5 \Delta u (x_1, x_2) - 0.5 \norm{\nabla u (x_1, x_2) }_2^2 ) &
 \text{if} \ x_2
 \ge (x_1-0.5)^2 + 0.1 \\
0.5 \Delta u (x_1, x_2) - 0.5 \norm{\nabla u (x_1, x_2) }_2^2 &
 \text{else} \enspace, \\
\end{array}
\right.
\] 
and for $(x_1, x_2) \in \partial \X$~: $\psi_1 (x_1, x_2) =
u(x_1,x_2)$ where 
\[ 
u (x_1, x_2) \, = \,\left\{
\begin{array}{l l}
(x_2 - ((x_1 - 0.5)^2+0.1))^3 &
 \text{if} \ x_2
 \ge (x_1-0.5)^2 + 0.1 \\
0 & \text{else} \enspace.
\end{array}
\right.
\] 
The definitions of the functions $f$, $\psi_1$ and $\psi_2$ are chosen
such that the function $u$, represented in Figure~\ref{solvipGraphSol},
is solution of~\eqref{VIex1} almost everywhere and such that the terms  
\nbcircle{1} and \nbcircle{2} in Equation~\eqref{VIex1} are non
positive for all $x \in \X$ (this condition must hold for the
variational inequality to be well-defined). 
This example leads to a free boundary problem for the actions of
\MAX. Indeed, the points of the state space $\X_h$ can be divided in
two parts, the points where \MAX\ chooses action $1$ (means
continue to play) and the points where \MAX\ chooses action $0$
(means that he stops the game). For $(x_1,x_2) \in \X_h$,  the optimal
strategy $\balpha$ for \MAX\ is $\balpha(x_1,x_2)=1$ if $x_2 \ge
(x_1-0.5)^2 + 0.1$ and $\balpha(x_1,x_2)=0$ else, for all $(x_1,x_2)\in\X$. 

As for the previous example, the domain $\X$ is discretized on a grid
with $m+1$ points in each direction, i.e.\ with a discretization step 
$h=\frac{1}{m}$ and we obtain a discrete space $\X_h$ with boundary $
\partial \X_h$. Then, Equation~\eqref{VIex1} is discretized by using the
discretization scheme~\eqref{decent}-~\eqref{centre}. 
After, the equations \nbcircle{1} and \nbcircle{2} are
simplified separately by keeping equations \eqref{visympli} true. In
this case, only equation \nbcircle{1} is multiply by $\frac{h^2}{c}$ with
$c$ an appropriate constant. After discretization, we obtain the
following dynamic programming equation for a game with state space $\X_h$~: 
\begin{align*}
v(x_i,x_j)  = \quad \max  \Bigg\{ & \min_{(b_1,b_2) \in \B} \
 \left( \frac{1}{2c} - b_1 \frac{h}{c}  \indicator_{b_1\le 0} \right) v(x_{i+1},x_j) 
 +  \left( \frac{1}{2c} + b_1 \frac{h}{c}  \indicator_{b_1 > 0} \right) v(x_{i-1},x_j) \\
 & +  \left( \frac{1}{2c} - b_2 \frac{h}{c}  \indicator_{b_2\le 0} \right) v(x_{i},x_{j+1}) 
 +  \left( \frac{1}{2c} + b_2 \frac{h}{c} \indicator_{b_2 > 0} \right) v(x_{i},x_{j-1}) \\
 & + \frac{h^2}{c} \,\frac{b_1^2 + b_2^2}{2} + \frac{h^2}{c} \,
 f(x_i,x_j),  \quad
 \psi_2 (x_i,x_j) \quad \Bigg\} \qquad \text{ for } (x_i,x_j) \in
 \X_h 
\end{align*}
with $c = 2 + h \abs{b_1} + h \abs{b_2} > 0$ and $v(x_i,x_j)
= \psi_1(x_i,x_j)$ for $(x_i,x_j)\in\partial\X_h $. 
The same comments about non-symmetry and the discount factor in
equation~\eqref{Iex1d} hold here. That is $\lambda =0 $ or equivalently
$\mu = 1$.

\begin{figure} 
\begin{tabular}{c c}
\includegraphics[height=4.5cm]{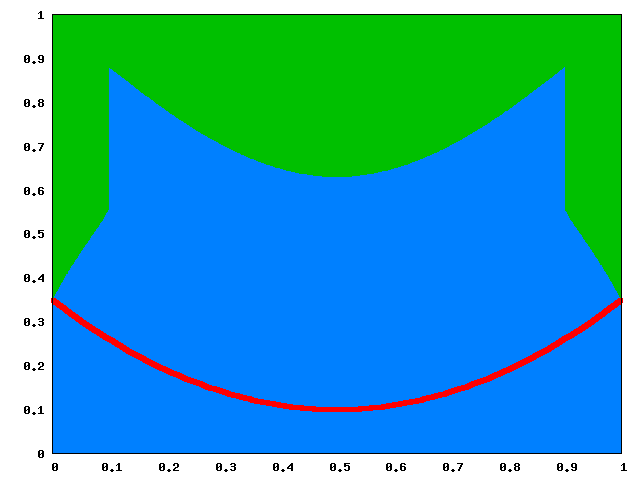} 
(a)
&
\includegraphics[height=4.5cm]{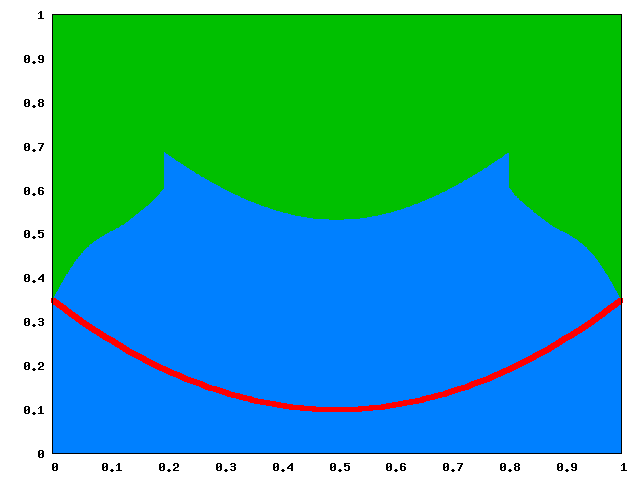} 
(b)
\\
\includegraphics[height=4.5cm]{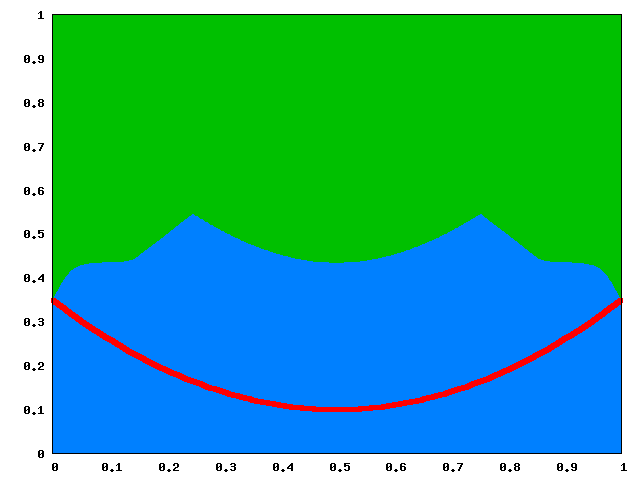} 
(c)
&
\includegraphics[height=4.5cm]{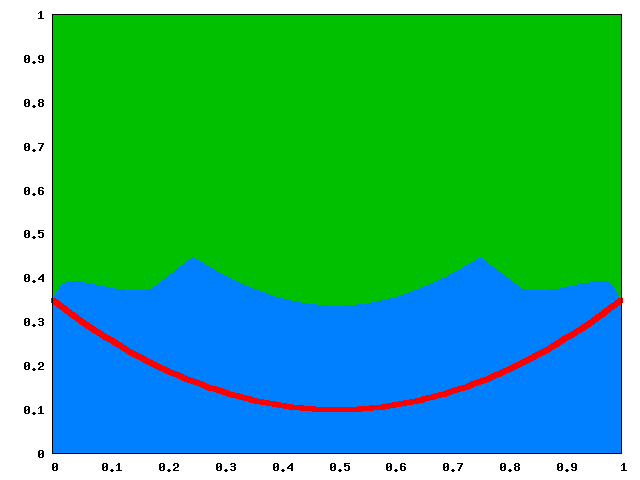} 
(d)
\\
\includegraphics[height=4.5cm]{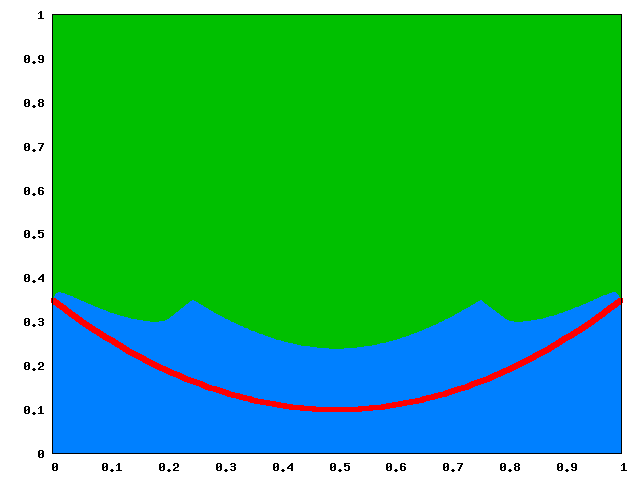} 
(e)
&
\includegraphics[height=4.5cm]{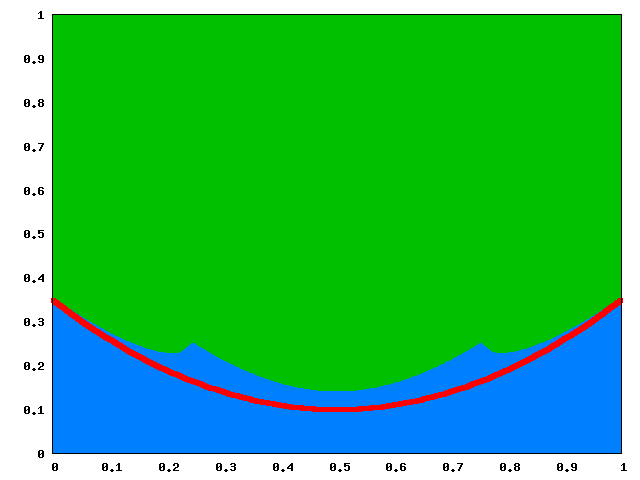} 
(f)
\\
\includegraphics[height=4.5cm]{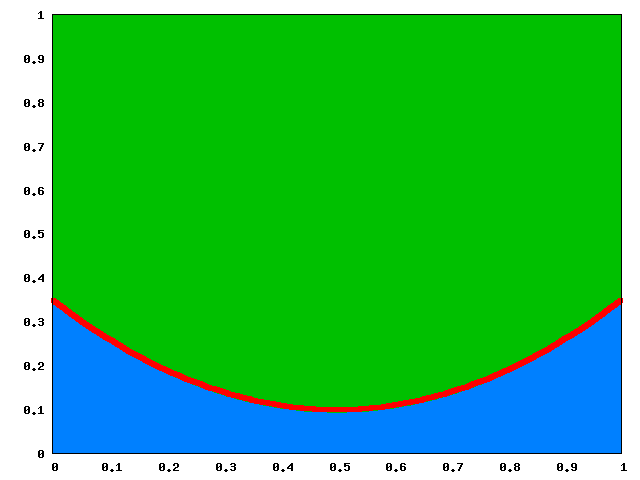} 
(g)
\end{tabular}
\caption{Application of AMG$\pi$ to the free boundary problem \eqref{VIex1} for a $1025 \times 1025$ points grid~: (a) after $100$ iterations, (b) after $200$ iterations, (c) after $300$ iterations, (d) after $400$ iterations, (e) after $500$ iterations, (f) after $600$ iterations and (g) after $700$ iterations.}
\label{solvipFaildTests}
\end{figure}

\begin{table} 
\centering
\caption{Numerical results for optimal stopping time game
  \eqref{VIex1} with a $1025 \times 1025$ points grid, computed by
  AMG$\pi$ with $\epsilon = 10^{-14}$.}
\label{VIex1t}
\begin{tabular}{| c| c| c| c| c| c| c|}
\hline 
 $\ki$ & $\nkj$ &  $\norm{r_v}_{\infty}$ &  $\norm{r_v}_{L_2}$ &  $\norm{e}_{\infty}$ & $\norm{e}_{L_2}$ & cpu time (s) \\
\hline 
$ 1 $ & $ 0 $ & $ 3.645e-01 $ & $ 9.195e-03 $ & $ 7.243e-01 $ & $ 1.998e-01 $ & $ 1.790e+00 $ \\
\hline 
$ 2 $ & $ 4 $ & $ 1.497e-01 $ & $ 1.347e-03 $ & $ 3.782e-01 $ & $ 1.218e-01 $ & $ 1.376e+01 $ \\
\hline 
$ 3 $ & $ 4 $ & $ 1.094e-01 $ & $ 8.839e-04 $ & $ 3.767e-01 $ & $ 1.213e-01 $ & $ 2.492e+01 $ \\
\hline 
\multicolumn{7}{|c|}{$\dots$} \\
\hline 
$ 100 $ & $ 3 $ & $ 1.744e-02 $ & $ 4.444e-05 $ & $ 2.392e-01 $ & $ 8.016e-02 $ & $ 1.009e+03 $ \\
\hline 
\multicolumn{7}{|c|}{$\dots$} \\
\hline 
$ 200 $ & $ 3 $ & $ 7.398e-03 $ & $ 1.879e-05 $ & $ 1.222e-01 $ & $ 3.996e-02 $ & $ 2.214e+03 $ \\
\hline 
\multicolumn{7}{|c|}{$\dots$} \\
\hline 
$ 300 $ & $ 3 $ & $ 2.510e-03 $ & $ 8.779e-06 $ & $ 5.614e-02 $ & $ 1.728e-02 $ & $ 3.619e+03 $ \\
\hline 
\multicolumn{7}{|c|}{$\dots$} \\
\hline 
$ 400 $ & $ 2 $ & $ 1.258e-03 $ & $ 4.363e-06 $ & $ 2.321e-02 $ & $ 6.519e-03 $ & $ 4.770e+03 $ \\
\hline 
\multicolumn{7}{|c|}{$\dots$} \\
\hline 
$ 500 $ & $ 2 $ & $ 4.761e-04 $ & $ 1.620e-06 $ & $ 6.601e-03 $ & $ 1.532e-03 $ & $ 5.861e+03 $ \\
\hline 
\multicolumn{7}{|c|}{$\dots$} \\
\hline 
$ 600 $ & $ 2 $ & $ 8.857e-05 $ & $ 2.781e-07 $ & $ 7.274e-04 $ & $ 9.598e-05 $ & $ 7.045e+03 $ \\
\hline 
\multicolumn{7}{|c|}{$\dots$} \\
\hline 
$ 650 $ & $ 2 $ & $ 1.533e-05 $ & $ 4.231e-08 $ & $ 1.538e-04 $ & $ 6.331e-05 $ & $ 7.630e+03 $ \\
\hline 
\multicolumn{7}{|c|}{$\dots$} \\
\hline 
$ 700 $ & $ 1 $ & $ 5.647e-08 $ & $ 8.734e-11 $ & $ 1.571e-04 $ & $ 6.619e-05 $ & $ 8.134e+03 $ \\
\hline 
$ 701 $ & $ 1 $ & $ 1.207e-08 $ & $ 2.267e-11 $ & $ 1.571e-04 $ & $ 6.619e-05 $ & $ 8.141e+03 $ \\
\hline 
$ 702 $ & $ 1 $ & $ 9.992e-16 $ & $ 7.284e-17 $ & $ 1.571e-04 $ & $ 6.619e-05 $ & $ 8.148e+03 $ \\
\hline 
\end{tabular}
\end{table}

\begin{table} 
\centering
\caption{Numerical results for optimal stopping time game
  \eqref{VIex1} with a $1025 \times 1025$ points grid, computed by
  FAMG$\pi$ with $c=10^{-2}$ and $\epsilon = 10^{-14}$.}
\label{VIex1FMG}
\begin{tabular}{| c| c| c| c| c| c| c|}
\hline
 $\ki$ & $\nkj$ &  $\norm{r_v}_{\infty}$ &  $\norm{r_v}_{L_2}$ &  $\norm{e}_{\infty}$ & $\norm{e}_{L_2}$ & cpu time (s) \\
\hline 
\multicolumn{7}{|c|}{ points in each direction : $ 3 $, step size : $ 5.00e-01 $ } \\ 
\hline 
$ 1 $ & $ 1 $ & $ 2.17e-01 $ & $ 2.17e-01 $ & $ 1.53e-01 $ & $ 1.53e-01 $ & $ << 1 $ \\
\hline 
$ 2 $ & $ 2 $ & $ 2.64e-05 $ & $ 2.64e-05 $ & $ 3.92e-02 $ & $ 3.92e-02 $ & $ << 1 $ \\
\hline 
\multicolumn{7}{|c|}{ points in each direction : $ 5 $, step size : $ 2.50e-01 $ } \\ 
\hline 
$ 1 $ & $ 2 $ & $ 2.19e-04 $ & $ 8.41e-05 $ & $ 3.02e-02 $ & $ 1.71e-02 $ & $ << 1 $ \\
\hline 
\multicolumn{7}{|c|}{ points in each direction : $ 9 $, step size : $ 1.25e-01 $ } \\ 
\hline 
$ 1 $ & $ 2 $ & $ 4.99e-03 $ & $ 1.06e-03 $ & $ 1.65e-02 $ & $ 7.99e-03 $ & $ << 1 $ \\
\hline 
$ 2 $ & $ 1 $ & $ 2.68e-03 $ & $ 5.41e-04 $ & $ 1.66e-02 $ & $ 8.15e-03 $ & $ << 1 $ \\
\hline 
$ 3 $ & $ 1 $ & $ 2.72e-04 $ & $ 5.49e-05 $ & $ 1.68e-02 $ & $ 8.30e-03 $ & $ << 1 $ \\
\hline 
\multicolumn{7}{|c|}{ points in each direction : $ 17 $, step size : $ 6.25e-02 $ } \\ 
\hline 
$ 1 $ & $ 2 $ & $ 2.26e-03 $ & $ 5.44e-04 $ & $ 8.75e-03 $ & $ 3.89e-03 $ & $ << 1 $ \\
\hline 
$ 2 $ & $ 1 $ & $ 7.97e-04 $ & $ 1.23e-04 $ & $ 8.84e-03 $ & $ 3.97e-03 $ & $ << 1 $ \\
\hline 
$ 3 $ & $ 1 $ & $ 4.65e-04 $ & $ 5.97e-05 $ & $ 8.98e-03 $ & $ 4.11e-03 $ & $ << 1 $ \\
\hline 
$ 4 $ & $ 1 $ & $ 9.57e-08 $ & $ 1.24e-08 $ & $ 9.01e-03 $ & $ 4.14e-03 $ & $ 1.00e-02 $ \\
\hline 
\multicolumn{7}{|c|}{ points in each direction : $ 33 $, step size : $ 3.12e-02 $ } \\ 
\hline 
$ 1 $ & $ 1 $ & $ 2.10e-04 $ & $ 1.90e-05 $ & $ 4.94e-03 $ & $ 2.16e-03 $ & $ 1.00e-02 $ \\
\hline 
$ 2 $ & $ 1 $ & $ 1.05e-04 $ & $ 6.57e-06 $ & $ 4.76e-03 $ & $ 2.09e-03 $ & $ 2.00e-02 $ \\
\hline 
\multicolumn{7}{|c|}{ points in each direction : $ 65 $, step size : $ 1.56e-02 $ } \\ 
\hline 
$ 1 $ & $ 1 $ & $ 6.26e-05 $ & $ 6.43e-06 $ & $ 2.49e-03 $ & $ 1.07e-03 $ & $ 4.00e-02 $ \\
\hline 
$ 2 $ & $ 1 $ & $ 3.64e-05 $ & $ 2.09e-06 $ & $ 2.45e-03 $ & $ 1.05e-03 $ & $ 7.00e-02 $ \\
\hline 
\multicolumn{7}{|c|}{ points in each direction : $ 129 $, step size : $ 7.81e-03 $ } \\ 
\hline 
$ 1 $ & $ 1 $ & $ 7.67e-06 $ & $ 3.88e-07 $ & $ 1.25e-03 $ & $ 5.33e-04 $ & $ 1.60e-01 $ \\
\hline 
\multicolumn{7}{|c|}{ points in each direction : $ 257 $, step size : $ 3.91e-03 $ } \\ 
\hline 
$ 1 $ & $ 1 $ & $ 2.86e-06 $ & $ 1.12e-07 $ & $ 6.28e-04 $ & $ 2.66e-04 $ & $ 6.20e-01 $ \\
\hline 
\multicolumn{7}{|c|}{ points in each direction : $ 513 $, step size : $ 1.95e-03 $ } \\ 
\hline 
$ 1 $ & $ 1 $ & $ 5.33e-07 $ & $ 1.44e-08 $ & $ 3.15e-04 $ & $ 1.33e-04 $ & $ 2.49e+00 $ \\
\hline 
\multicolumn{7}{|c|}{ points in each direction : $ 1025 $, step size : $ 9.77e-04 $ } \\ 
\hline 
$ 1 $ & $ 2 $ & $ 1.79e-07 $ & $ 3.82e-09 $ & $ 1.57e-04 $ & $ 6.62e-05 $ & $ 1.58e+01 $ \\
\hline 
$ 2 $ & $ 1 $ & $ 9.66e-08 $ & $ 8.84e-10 $ & $ 1.57e-04 $ & $ 6.62e-05 $ & $ 2.30e+01 $ \\
\hline 
$ 3 $ & $ 1 $ & $ 5.39e-08 $ & $ 4.10e-10 $ & $ 1.57e-04 $ & $ 6.62e-05 $ & $ 3.00e+01 $ \\
\hline 
$ 4 $ & $ 1 $ & $ 2.86e-08 $ & $ 1.31e-10 $ & $ 1.57e-04 $ & $ 6.62e-05 $ & $ 3.70e+01 $ \\
\hline 
$ 5 $ & $ 1 $ & $ 7.41e-09 $ & $ 1.60e-11 $ & $ 1.57e-04 $ & $ 6.62e-05 $ & $ 4.34e+01 $ \\
\hline 
$ 6 $ & $ 1 $ & $ 8.88e-16 $ & $ 7.31e-17 $ & $ 1.57e-04 $ & $ 6.62e-05 $ & $ 4.99e+01 $ \\
\hline 
\end{tabular}
\end{table}

\begin{figure} 
\includegraphics[height=4.5cm]{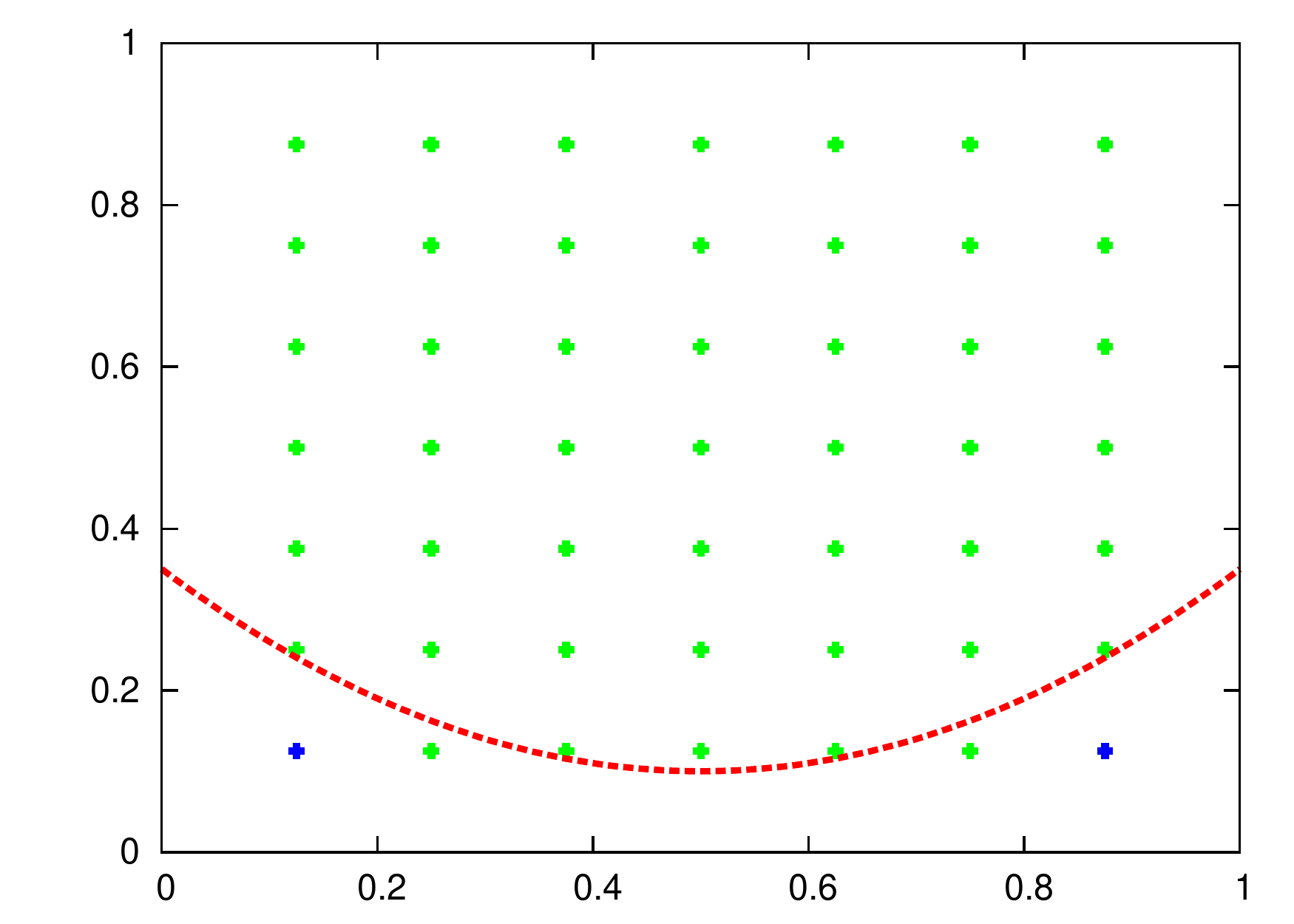}(a)
\includegraphics[height=4.5cm]{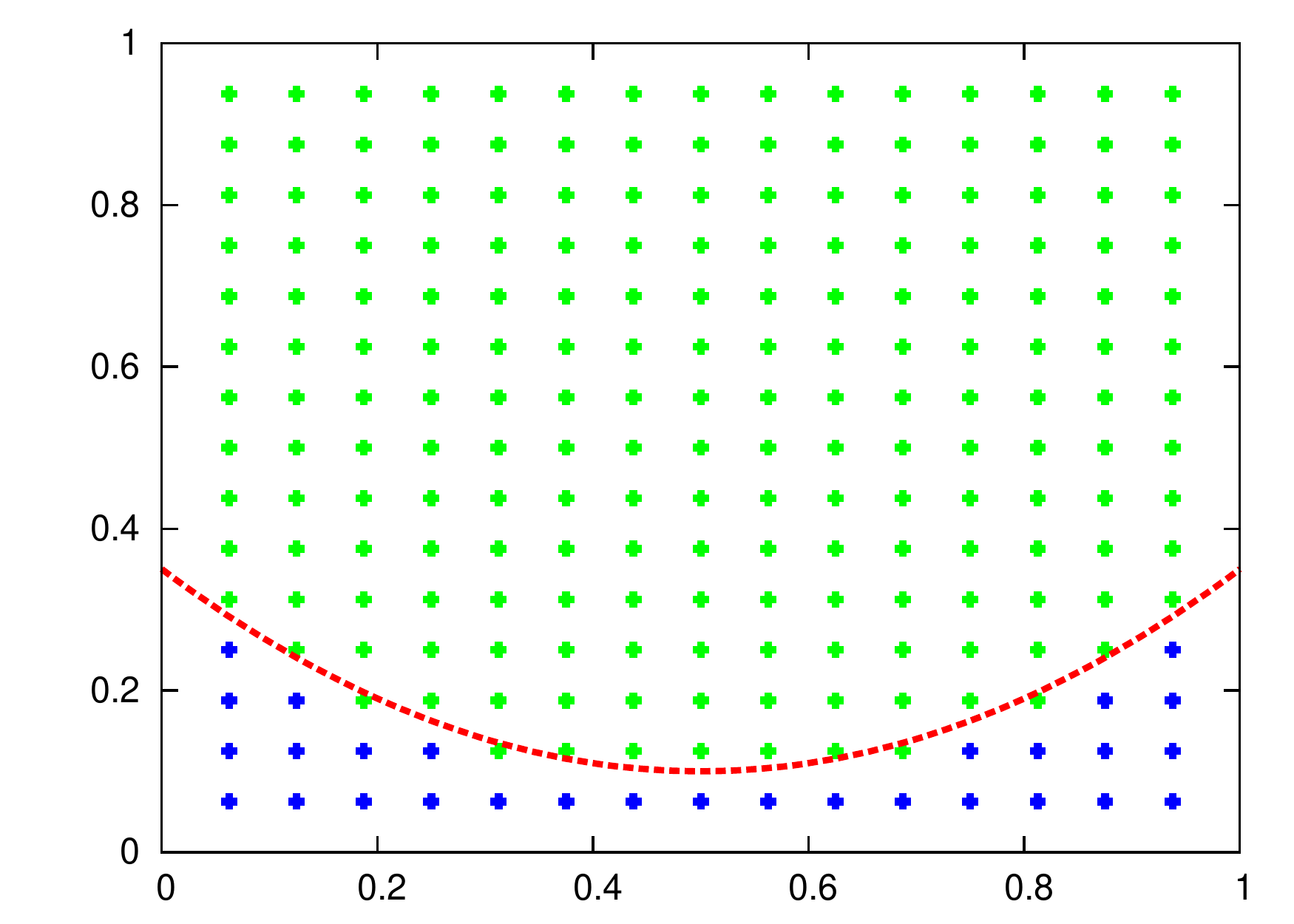}(b)

\includegraphics[height=4.5cm]{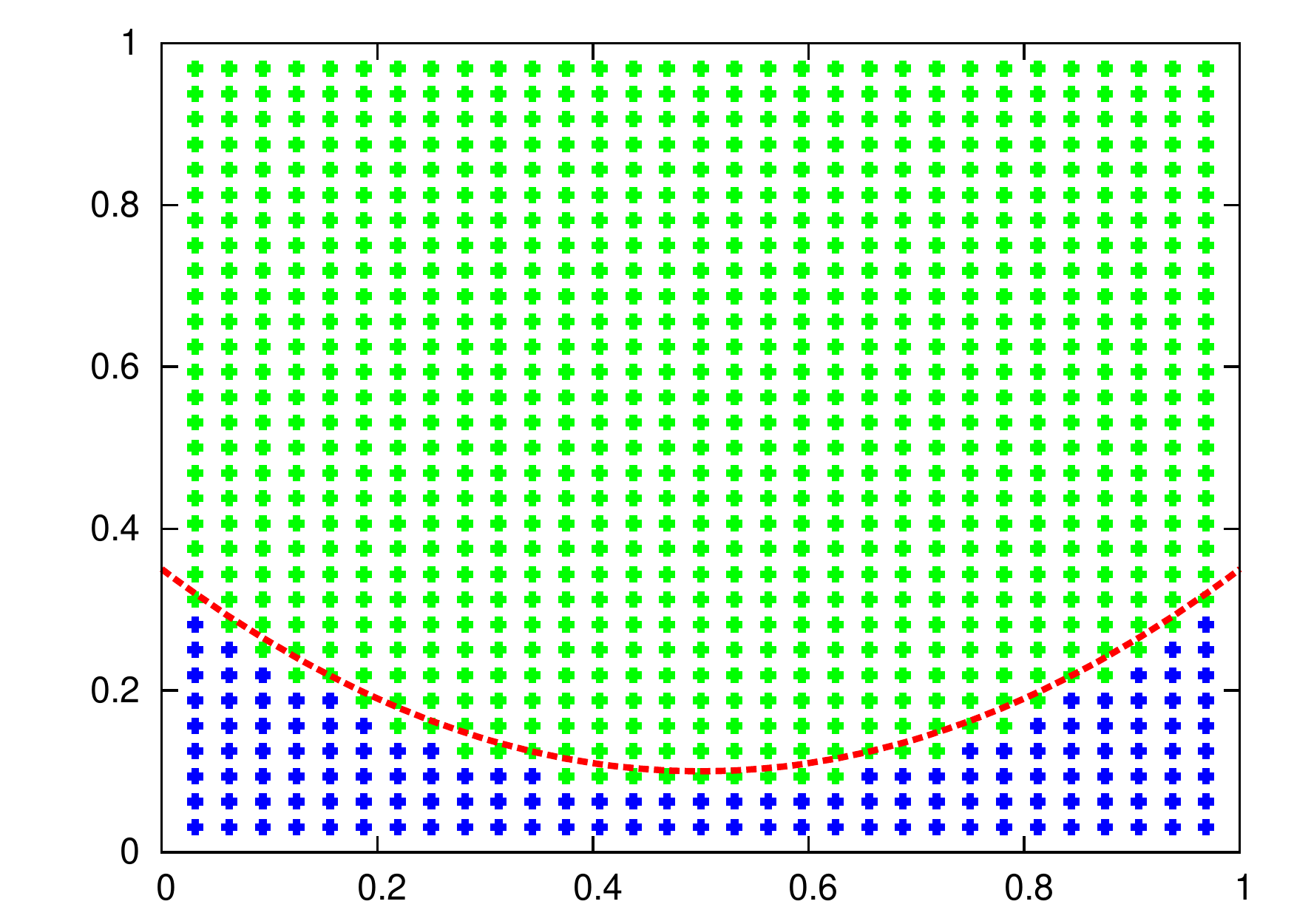}(c)
\includegraphics[height=4.5cm]{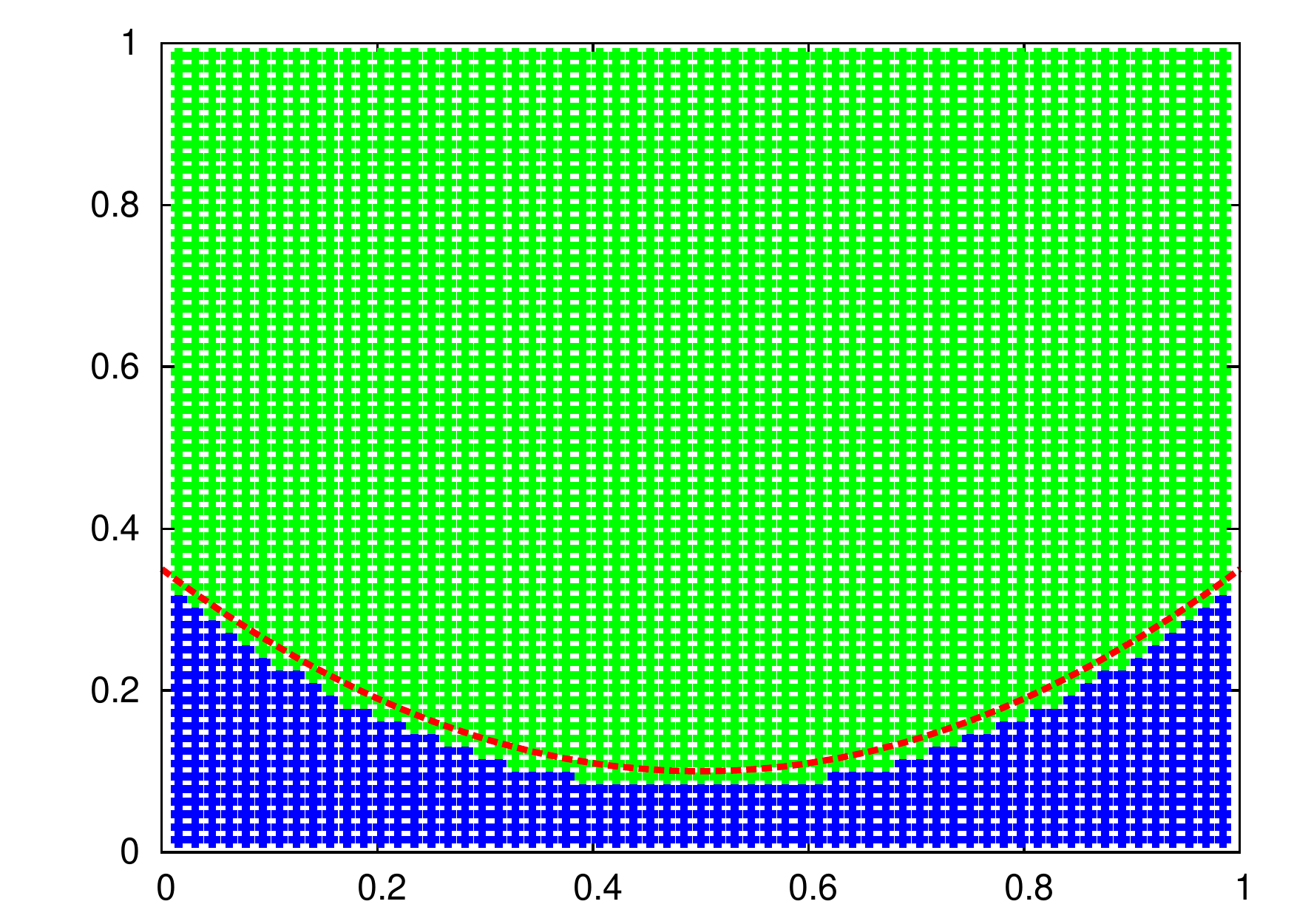}(d)
\caption{Application FAMG$\pi$ to the free boundary problem \eqref{VIex1} for: (a) $9 \times 9$ points grid, (b) $17 \times 17$ points grid, (c) $33 \times 33$ points grid, (d) $65 \times 65$ points grid.}
\label{solvipSuccessTests}
\end{figure}

The numerical results are performed for Equation \eqref{VIex1} when
discretized on a grid with $1025$ points in each direction. In the
domain $\X_h$, for a fixed strategy $\balpha$ of \MAX, we represent a
points $x$ with a green color when 
$\balpha(x) = 1$, that is where \MAX\ decides to continue
playing, and with a blue color when $\balpha(x) = 0$, that is when
\MAX\ decides to stop the game. 
The optimal strategy for \MAX\ is to have only
green points above the red curve, $x_2= (x_1-0.5)^2 + 0.1$, and only
blue points under.
We start the tests with $\balpha(x) = 0$ for all $x \in\X$, that is
with blue points in the whole domain. 

Numerical results with AMG$\pi$ are shown geometrically in
Figure~\ref{solvipFaildTests} where the strategies of \MAX\
obtained after $100$, $200$, $300$, $400$, $500$, $600$ and $700$
iterations are represented. 
We observe in Table~\ref{VIex1t} that AMG$\pi$ finds an approximation of the
solution after $702$ iterations and in about two hours and
15 minutes. The stopping criterion for policy iterations of
AMG$\pi$ in this test is $\epsilon = 10^{14}$. This criterion was
chosen to ensure the convergence of the policy iterations, indeed with
a smaller $\epsilon$ it did not converge because the intern policy
iterations did not gave a precise enough approximation. 

In table \ref{VIex1FMG}, we present numerical results for the
application of  FAMG$\pi$ with $c=10^{-2}$ and $\epsilon = 10^{-14}$
to problem \eqref{VIex1} for a $1025 \times 1025$ points grid.  
We observe that our algorithm solves the problem in about $49$ seconds. 
Geometrical representation of the strategies of \MAX\, obtained by AMG$\pi$ 
on four successive levels in the FAMG$\pi$ algorithm, 
are shown in Figure \ref{solvipSuccessTests}.  
We can see that on coarse grids, the algorithm can find a good
approximation of the solution in a few iterations.  
The interpolation of this solution and the corresponding strategies,
are used to start AMG$\pi$ on the next fine level and we observe that only a few
numbers of policy iterations are needed on each level. 

With this example we show the advantage of using FAMG$\pi$.  
Indeed, the computation time of the FAMG$\pi$ algorithm
seems to be in the order of the number of discretization points whereas
that of a AMG$\pi$ algorithm is about $160$ times greater. 
This is due to the large number of iterations needed by AMG$\pi$ for
solving this kind of games. 
Indeed, this number should be compared to the diameter of the graph 
(that is the largest number of edges which must be cover to travel
from one point to another) associated to the corresponding game
problem, for instance the union of all graphs of the Markov 
chains associated to all couple of fixed policies $\balpha$ and $\bbeta$.
Hence due to the finite differences discretization, the arcs of the
graphs are supported by edges of the grids $\X_h$ in $\Z^2$, so the
diameter is $2m$ with $m = 1024$.

\subsection{Stopping game with two optimal stopping}

In this example, we consider a stopping game where both players have
the possibility to stop the game, see~\cite{FriedmanA} for a complete theory about
this subject. In this case, the value of the game
starting in $x\in\X$ is given by~:
\[ 
  v(x) \, = \, \sup_{\stops_1} \ \inf_{\stops_2} \, \Big\{ \,\sE^{\stops_1, \stops_2}_{x} \left[\, 
 \int^{\stops_1 \wedge \stops_2}_{0} r(\Xk_t, \Bk_t)\,dt +
 \psi_1(\Xk_{\stops_1}) \, \indicator_{\stops_1 < \stops_2} +
 \psi_2(\Xk_{\stops_2}) \, \indicator_{\stops_2 \le \stops_1} \, \Big|
 \, \Xk_0 =x \,\right] \, \Big\} 
\] 
where $\stops_1 \wedge \stops_2 = \min{(\stops_1, \stops_2)}$ and we
assume $\min{(\stops_1, \stops_2)} < \tau$ ($\tau = \inf \set{t \ge 0|\Xk_t \notin \X}$, 
then $v$ is 
solution of equation~:
\begin{equation}\label{VI2stop}
\max \ \bigg\{ \, \psi_1 (x)  - v (x), \,
\min \left\{ \psi_2 (x)  - v (x)\,, \,  L(v; x) + r(x) \right\}
\, \bigg\} \,=\, 0 \qquad \text{for $x$ in $\X,$}
\end{equation}
or equivalently, 
\[
  \left\{
  \begin{array}{l l}
    (L(v; x) + r(x)) (w(x) - v(x)) \le 0 &\text{for $x \in \X$},\\
    \forall w,  \ \psi_1 \le w \le \psi_2
    \text{ and } \psi_1 \le v \le \psi_2 &  , \\
  \end{array}
  \right.
\]
that is 
\[
\text{for $x \in \X$} \quad
  \left\{
  \begin{array}{l l}
    (L(v; x) + r(x)) \le 0 & \text{ if } v(x) = \psi_1(x) \\
    (L(v; x) + r(x)) \ge 0 & \text{ if } v(x) = \psi_2(x)\\
    (L(v; x) + r(x)) = 0 & \text{ if } \psi_1(x) < v(x) < \psi_2(x).\\
  \end{array}
  \right.
\]

\begin{figure}
\centering
\scalebox{0.8}{
\input{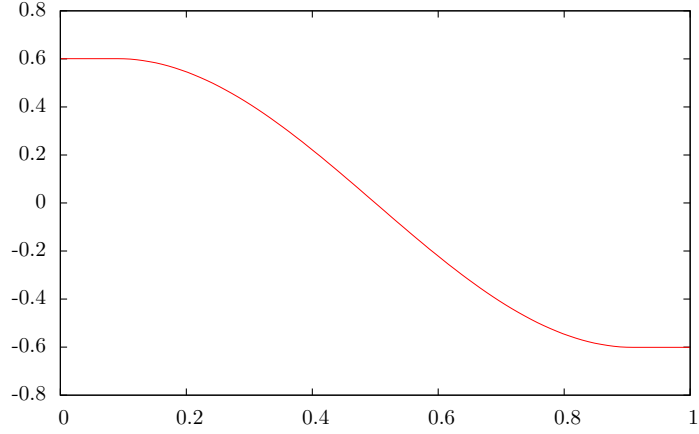}
}
\caption{Solution of Equation~\eqref{exVI2stop}}
\label{fig-exVI2stop-sol}
\end{figure}

For the numerical tests, we consider the stochastic differential game
whose value $v$ is solution of~:
\begin{equation}\label{exVI2stop}
\max \ \bigg\{ \, \psi_1 (x)  - v (x), \,
\min \left\{ \psi_2 (x)  - v (x)\,, \,  0.5 \, \nabla(x) + r(x) \right\}
\, \bigg\} \,=\, 0 \qquad \text{for $x$ in $\X,$}
\end{equation}
where $\X = [0,1]$, 
for all $x \in \X$: $\psi_1(x)=-\bar \psi_2$, $\psi_2(x)=\bar \psi_2$
  with $\bar \psi_2 = (2 \cos(0.09\pi) + \pi
 (0.18 -1) \sin(0.09\pi)) / 2) \approx 0.6$ and 
$r (x) \, = \, 0.5 \, \pi^2 \, \cos(\pi x)$. 
For all $x\in\X$, the sets of actions are $\A = \set{0, 1}$ for \MAX\
and  $\B = \set{0,  1}$ for \MIN, where action $0$ means that the
player chooses to stop the game and receive $\psi_1$ when \MAX\
stops or $\psi_2$ when \MIN\ stops, action $1$ means that the game is
continuing.   
Here, the exact solution of Equation~\eqref{exVI2stop} in the viscosity sense is
\[
\text{for $x \in \X$} \quad
  \left\{
  \begin{array}{l l}
    \psi_1(x) & \text{ for } x > (1 - 0.09)\\
    \psi_2(x) & \text{ for } x < 0.09 \\
    \cos(\pi x) +  \pi \sin(0.09\pi ) x + c & \text{ for } 0.09 > x > (1 - 0.09) \\
  \end{array}
  \right.
\]
where the constant $c = (\bar \psi_2 - \cos(0.09\pi) - 0.09 \pi\sin(0.09\pi))$ and
is represented in Figure~\ref{fig-exVI2stop-sol}. 
For all $x\in\X$, the optimal strategy for \MAX\ is $\balpha(x) = 0$
if $x > (1 - 0.09)$ and $\balpha(x) = 1$ else. For all $x\in\X$, the
optimal strategy for \MIN\ is $\bbeta(x) = 0$ if $x < (0.09)$ and
$\bbeta(x) = 1$ else.  

We present numerical results for the discretization of
Equation~\eqref{exVI2stop} on a grid with $2049$ points 
in Table~\ref{VIex2AMGpi-tab} when using
AMG$\pi$ with $\epsilon = 10^{-10}$ and in Table~\ref{VIex2FMG-tab}
when using FAMG$\pi$ with $c=10^{-2}$ and $\epsilon = 10^{-10}$.
As in the previous example, we see the advantage of using FAMG$\pi$
for this kind of games.
Indeed, FAMG$\pi$ solves the problem in about one second  
while AMG$\pi$ needs about $24$ minutes. 
As for the previous example, the computation time of the FAMG$\pi$ seems to
be in the order of the number of discretization points.
For this example, due to the finite differences
discretization, the diameter of the graph is $m$ with $m = 2048$. We
see in Table~\ref{VIex2AMGpi-tab} that both numbers of intern and
external policy iterations for AMG$\pi$ are of the order of the diameter of the graph.

\begin{table}
\centering
\caption{Numerical results for optimal stopping time game
  \eqref{exVI2stop} with a $2049 \times 2049$ points grid, computed by
  AMG$\pi$ with $\epsilon = 10^{-10}$.}
\label{VIex2AMGpi-tab}
\begin{tabular}{| c | c| c| c| c| c| c| c|} 
\hline 
 $\ki$ & $\nkj$ &  $\norm{r_v}_{\infty}$ &  $\norm{r_v}_{L_2}$ &  $\norm{e}_{\infty}$ & $\norm{e}_{L_2}$ & cpu time (s) \\ 
\hline 
$1$  & $1$  & $1.20e+00$  & $7.80e-01$  & $7.75e-01$  & $2.91e-01$  & $<< 1$   \\ 
 \hline 
$2$  & $863$  & $1.20e+00$  & $2.66e-02$  & $6.02e-01$  & $3.00e-01$  & $1.58e+00$   \\ 
 \hline 
$3$  & $1025$  & $1.20e+00$  & $2.66e-02$  & $6.03e-01$  & $3.00e-01$  & $3.31e+00$   \\ 
 \hline 
\multicolumn{7}{|c|}{$\dots$} \\ 
 \hline 
$100$  & $1026$  & $1.20e+00$  & $2.66e-02$  & $7.10e-01$  & $2.46e-01$  & $1.81e+02$   \\ 
 \hline 
\multicolumn{7}{|c|}{$\dots$} \\ 
 \hline 
$200$  & $992$  & $1.20e+00$  & $2.66e-02$  & $8.16e-01$  & $1.93e-01$  & $3.75e+02$   \\ 
 \hline 
\multicolumn{7}{|c|}{$\dots$} \\ 
 \hline 
$300$  & $947$  & $1.20e+00$  & $2.66e-02$  & $9.16e-01$  & $1.44e-01$  & $5.72e+02$   \\ 
 \hline 
\multicolumn{7}{|c|}{$\dots$} \\ 
 \hline 
$400$  & $910$  & $1.20e+00$  & $2.66e-02$  & $1.00e+00$  & $1.01e-01$  & $7.73e+02$   \\ 
 \hline 
\multicolumn{7}{|c|}{$\dots$} \\ 
 \hline 
$500$  & $882$  & $1.20e+00$  & $2.66e-02$  & $1.08e+00$  & $6.59e-02$  & $9.78e+02$   \\ 
 \hline 
\multicolumn{7}{|c|}{$\dots$} \\ 
 \hline 
$600$  & $862$  & $1.20e+00$  & $2.66e-02$  & $1.14e+00$  & $4.05e-02$  & $1.19e+03$   \\ 
 \hline 
\multicolumn{7}{|c|}{$\dots$} \\ 
 \hline 
$700$  & $849$  & $1.20e+00$  & $2.66e-02$  & $1.18e+00$  & $2.84e-02$  & $1.41e+03$   \\ 
 \hline 
\multicolumn{7}{|c|}{$\dots$} \\ 
 \hline 
$800$  & $843$  & $1.20e+00$  & $2.66e-02$  & $1.20e+00$  & $2.65e-02$  & $1.64e+03$   \\ 
 \hline 
\multicolumn{7}{|c|}{$\dots$} \\ 
 \hline 
$839$  & $843$  & $1.20e+00$  & $2.66e-02$  & $1.20e+00$  & $2.66e-02$  & $1.73e+03$   \\ 
 \hline 
$840$  & $843$  & $2.03e-07$  & $4.50e-09$  & $5.22e-07$  & $2.57e-07$  & $1.74e+03$   \\ 
 \hline 
$841$  & $1$  & $1.11e-16$  & $6.57e-18$  & $1.16e-07$  & $7.40e-08$  & $1.74e+03$   \\ 
 \hline 
\end{tabular} 
\end{table}

\begin{table}
\centering
\caption{Numerical results for optimal stopping time game
  \eqref{exVI2stop} with a $2049 \times 2049$ points grid, computed by
  FAMG$\pi$ with $c=10^{-2}$ and $\epsilon = 10^{-10}$.}
\label{VIex2FMG-tab}
\begin{tabular}{| c | c| c| c| c| c| c| c|} 
\hline 
 $\ki$ & $\nkj$ &  $\norm{r_v}_{\infty}$ &  $\norm{r_v}_{L_2}$ &  $\norm{e}_{\infty}$ & $\norm{e}_{L_2}$ & cpu time (s) \\ 
\hline 
\multicolumn{7}{|c|}{ points in each direction : $ 3 $, step size : $ 5.00e-01 $ } \\ 
 \hline 
$1$  & $2$  & $1.20e+00$  & $1.20e+00$  & $6.01e-01$  & $6.01e-01$  & $<< 1$   \\ 
 \hline 
$2$  & $2$  & $0.00e+00$  & $0.00e+00$  & $5.55e-17$  & $5.55e-17$  & $<< 1$   \\ 
 \hline 
\multicolumn{7}{|c|}{ points in each direction : $ 5 $, step size : $ 2.50e-01 $ } \\ 
 \hline 
$1$  & $2$  & $1.56e-01$  & $9.02e-02$  & $1.13e-01$  & $1.03e-01$  & $<< 1$   \\ 
 \hline 
$2$  & $1$  & $1.89e-17$  & $1.09e-17$  & $1.13e-01$  & $9.22e-02$  & $<< 1$   \\ 
 \hline 
\multicolumn{7}{|c|}{ points in each direction : $ 9 $, step size : $ 1.25e-01 $ } \\ 
 \hline 
$1$  & $2$  & $1.20e+00$  & $4.54e-01$  & $8.74e-01$  & $3.52e-01$  & $<< 1$   \\ 
 \hline 
$2$  & $5$  & $1.20e+00$  & $4.54e-01$  & $1.09e+00$  & $4.14e-01$  & $<< 1$   \\ 
 \hline 
$3$  & $5$  & $5.55e-17$  & $2.21e-17$  & $5.74e-03$  & $4.35e-03$  & $<< 1$   \\ 
 \hline 
\multicolumn{7}{|c|}{ points in each direction : $ 17 $, step size : $ 6.25e-02 $ } \\ 
 \hline 
$1$  & $2$  & $1.20e+00$  & $3.10e-01$  & $1.16e+00$  & $3.00e-01$  & $<< 1$   \\ 
 \hline 
$2$  & $10$  & $1.28e-03$  & $3.30e-04$  & $5.19e-03$  & $2.83e-03$  & $<< 1$   \\ 
 \hline 
$3$  & $1$  & $0.00e+00$  & $0.00e+00$  & $3.36e-03$  & $2.10e-03$  & $<< 1$   \\ 
 \hline 
\multicolumn{7}{|c|}{ points in each direction : $ 33 $, step size : $ 3.12e-02 $ } \\ 
 \hline 
$1$  & $2$  & $0.00e+00$  & $0.00e+00$  & $1.36e-04$  & $9.50e-05$  & $<< 1$   \\ 
 \hline 
\multicolumn{7}{|c|}{ points in each direction : $ 65 $, step size : $ 1.56e-02 $ } \\ 
 \hline 
$1$  & $2$  & $1.20e+00$  & $1.51e-01$  & $1.20e+00$  & $1.51e-01$  & $<< 1$   \\ 
 \hline 
$2$  & $28$  & $0.00e+00$  & $0.00e+00$  & $7.08e-05$  & $4.94e-05$  & $<< 1$   \\ 
 \hline 
\multicolumn{7}{|c|}{ points in each direction : $ 129 $, step size : $ 7.81e-03 $ } \\ 
 \hline 
$1$  & $2$  & $1.20e+00$  & $1.07e-01$  & $1.20e+00$  & $1.07e-01$  & $<< 1$   \\ 
 \hline 
$2$  & $54$  & $2.78e-17$  & $3.75e-18$  & $6.66e-05$  & $3.85e-05$  & $2.00e-02$   \\ 
 \hline 
\multicolumn{7}{|c|}{ points in each direction : $ 257 $, step size : $ 3.91e-03 $ } \\ 
 \hline 
$1$  & $2$  & $1.20e+00$  & $7.53e-02$  & $1.20e+00$  & $7.52e-02$  & $2.00e-02$   \\ 
 \hline 
$2$  & $108$  & $1.20e+00$  & $7.53e-02$  & $1.20e+00$  & $7.53e-02$  & $7.00e-02$   \\ 
 \hline 
$3$  & $107$  & $1.11e-16$  & $8.53e-18$  & $1.61e-06$  & $1.05e-06$  & $1.30e-01$   \\ 
 \hline 
\multicolumn{7}{|c|}{ points in each direction : $ 513 $, step size : $ 1.95e-03 $ } \\ 
 \hline 
$1$  & $2$  & $1.20e+00$  & $5.32e-02$  & $1.20e+00$  & $5.32e-02$  & $1.30e-01$   \\ 
 \hline 
$2$  & $212$  & $1.11e-16$  & $9.82e-18$  & $4.53e-07$  & $3.02e-07$  & $3.00e-01$   \\ 
 \hline 
\multicolumn{7}{|c|}{ points in each direction : $ 1025 $, step size : $ 9.77e-04 $ } \\ 
 \hline 
$1$  & $2$  & $1.20e+00$  & $3.76e-02$  & $1.20e+00$  & $3.76e-02$  & $3.00e-01$   \\ 
 \hline 
$2$  & $422$  & $1.11e-16$  & $1.15e-17$  & $1.69e-07$  & $1.18e-07$  & $9.40e-01$   \\ 
 \hline 
\multicolumn{7}{|c|}{ points in each direction : $ 2049 $, step size : $ 4.88e-04 $ } \\ 
 \hline 
$1$  & $2$  & $2.03e-07$  & $4.50e-09$  & $5.22e-07$  & $2.57e-07$  & $9.40e-01$   \\ 
 \hline 
$2$  & $1$  & $1.11e-16$  & $8.52e-18$  & $1.16e-07$  & $7.40e-08$  & $9.50e-01$   \\ 
 \hline 
\end{tabular} 
\end{table}

\section{Conclusion and perspective}

In this paper, we have presented our algorithm AMG$\pi$ for solving two
player zero-sum stochastic games. This program combines the policy
iteration algorithm with algebraic multigrid methods. 
Our experiences on a Isaacs equation show better results for AMG$\pi$ in comparison
with policy iteration combined to a direct linear solver. 
We observed that the most part of the computation time for the
resolution of a non-linear equation~\eqref{eq1} is used to solved
the linear systems involved in the policy iteration algorithm. 
Hence, we noticed that the computation time of AMG$\pi$ increase linearly
with the size of the problem.  

Furthermore, we also presented a full multi-level algorithm, called FAMG$\pi$,
for solving two player zero-sum stochastic differential games.  
The numerical results on some stopping differential stochastic games
presented here show that FAMG$\pi$ improves substantially the
computation time of the policy iteration algorithm for this kind of games.  
Indeed the computation time of FAMG$\pi$ 
seems to be in the order of the number of discretization points whereas
that of AMG$\pi$ algorithm is about $160$ to $1700$ times greater.
This is due to the large number of iterations needed by AMG$\pi$ for
solving this kind of games. Indeed, this number
should be compared to the diameter of the graph associated to the
corresponding game problem, for instance the union of all graphs of
the Markov chains associated to fixed policies $\balpha$ and $\bbeta$.

The FAMG$\pi$ algorithm uses coarse grids discretizations of
the partial differential equation and so cannot be applied directly to
the dynamic programming equation of a two player zero-sum stochastic
game with finite state space.  
One may ask if adapting the FAMG$\pi$ algorithm to
this kind of games is possible.  
Indeed, the complexity of two player zero-sum stochastic games is
still unsettled, one only knows that it belongs to the complexity
class of NP$\cap$coNP~\cite{Puri}, and any new approach maybe useful
to understand this complexity. 



\linespread{1}

\bibliographystyle{plain}

\bibliography{references}

\end{document}